\documentclass[11pt]{asaproc} 
\usepackage{amsmath}
\usepackage{amsfonts}
\usepackage{amssymb}
\usepackage{epsfig}
\usepackage{makeidx}
\usepackage{natbib}
\usepackage{color}
\usepackage{amsthm}
\usepackage{framed}
\usepackage{float}
\usepackage{booktabs}
\usepackage{siunitx}
\usepackage{graphicx}
\usepackage{mathtools}
\usepackage{multirow, makecell}
\usepackage[normalem]{ulem}

\usepackage{geometry}
\usepackage{times}

\bibpunct{(}{)}{;}{a}{,}{,}

\usepackage{color}

\newcommand{\bn}{\begin{enumerate}}
\newcommand{\en}{\end{enumerate}}
\newcommand{\bi}{\begin{itemize}}
\newcommand{\ei}{\end{itemize}}
\newcommand{\be}{\begin{eqnarray}}
\newcommand{\ee}{\end{eqnarray}}
\newcommand{\by}{\begin{eqnarray*}}
\newcommand{\ey}{\end{eqnarray*}}
\newcommand{\beq}{\begin{equation}}
\newcommand{\eeq}{\end{equation}}

 \title{Exponentiated Generalized Pareto Distribution: Properties and applications towards Extreme Value Theory}
\author{Se Yoon Lee \thanks{Texas A \& M University, 400 Bizzell St, College Station, TX 77843} \and  Joseph H. T. Kim \thanks{Yonsei University, Seoul, Korea, 120-749}}

\begin{document}
\maketitle
\newtheorem{thm}{Theorem}
\newtheorem{cor}[thm]{Corollary}
\newtheorem{lem}[thm]{Lemma}
\newtheorem{prop}[thm]{Proposition}
\newtheorem{defn}[thm]{Definition}
\newtheorem{exam}[thm]{Example}

\begin{abstract}
\noindent The Generalized Pareto Distribution (GPD) plays a central role in modelling heavy tail phenomena in many applications. Applying the GPD to actual datasets however is a non-trivial task. One common way suggested in the literature to investigate the tail behaviour is to take logarithm to the original dataset in order to reduce the sample variability. Inspired by this, we propose  and study the Exponentiated Generalized Pareto Distribution (exGPD), which is created via log-transform of the GPD variable. After introducing the exGPD we derive various distributional quantities, including the moment generating function, tail risk measures.  As an application we also develop a plot as an alternative to the Hill plot to identify the tail index of heavy tailed datasets, based on the moment matching for the exGPD. Various numerical analyses with both simulated and actual datasets show that the proposed plot works well.
\end{abstract}

\section{Introduction}
The Generalized Pareto Distribution (GPD) has recently emerged as an influential distribution in modelling heavy tailed datasets in various applications in finance, operational risk, insurance and environmental studies. In particular, it is widely used to model sample exceedances beyond some large threshold, a procedure commonly known as the peaks-over-threshold (POT) method in the Extreme Value Theory (EVT) literature;  see, e.g.,  \cite{EmbrechtsModellingExtremalEvents} and \cite{Beirl+Goege+etal:06a}. The heavy tail phenomenon and the GPD are linked through the famous Pickands-Balkema-de Hann theorem (\cite{balkema1974residual} and \cite{pickands1975statistical}) which states that, for an arbitrary distribution of which the sample maximum tends to a non-degenerate distribution after suitable standardization, the distribution function of its exceedances over a large threshold converges to the GPD. To this extent, extensive research has been carried out in the literature to characterize the GPD and apply it to the EVT framework; see, for example, \cite{de2010parameter} for a comprehensive survey on various estimation procedures for the GPD parameter.\\

In this paper we introduce a two-parameter exponentiated Generalized Pareto Distribution (exGPD in short) and study its distributional properties, which is the first contribution of the present paper. The term `exponentiated' is used because this new distribution is obtained via log-transform of the GPD random variable. There are ample examples where a new distribution is constructed through logarithm or exponential transformation of an existing one. Such examples include normal and log-normal, gamma and log-gamma, and Pareto and shifted exponential distributions. Thus introducing the exGPD is interesting on its own right from statistical viewpoint, but we have further motivations of considering such a distribution in connection to EVT. In various graphical tools offered by EVT one frequently investigates the tail behaviour with log-transformed data rather than the original data, as seen in, for example, \cite{Hill:75a} plot and the estimator of  \cite{pickands1975statistical}. As log transformation greatly reduces the variability of extreme quantiles, this practice naturally allows one to investigate the tail behaviour in a more stable manner. Therefore, given that the GPD is the central distribution in EVT modelling, it makes sense to create an alternative plotting tool using the exGPD directly. This is our second contribution to this paper. In particular, we develop a new plot as an alternative to the Hill plot to identify the tail index of heavy tailed datasets. The proposed Log Variance (LV) plot is based on the idea of the sample variance of log exceedances to be matched to the variance of the exGPD. Through various numerical illustrations with both simulated and actual datasets, it is shown that the LV plot works reasonably well compared to the Hill plot, elucidating the usefulness of the exGPD. \\

This article is organized as follows. In Section 2 we define the exGPD and investigate its distributional quantities including the moment generating function, from which the moments can be obtained. We show that the moment of all orders are finite for the exGPD, unlike the GPD. Section 3 derives popular risk measures, including the Value-at-Risk and the conditional tail expectation. In Section 4 we develop  the LV plot that is derived from the exGPD variance. We use  both simulated and actual datasets to illustrate the proposed plot and compare it to the Hill plot. Section 5 concludes the article.

\section{Exponentiated Generalized Pareto Distribution}
\subsection{Definition}
The distribution function (df) of the two-parameter GPD with parameter ($\sigma, \xi$) is defined as 
\begin{equation}\label{def.GPD.df1}
G_{X}(x) = 1 - \Big(1 + \dfrac{\xi x}{\sigma} \Big)^{-1/\xi} ; \quad \xi \neq 0
\end{equation}
where the support is $x \geq 0$ for $\xi > 0$ and $0 \leq x \leq -\sigma/\xi$ for $\xi < 0$. Here, $\sigma$ and $\xi$ are called the scale and shape parameter, respectively. For the case of $\xi = 0$, the df is defined as
\begin{equation}
G_{X}(x) = 1 - e^{-x/\sigma}; \quad \xi=0,
\end{equation} an exponential distribution defined on $ x \geq 0$ with scale parameter $\sigma$. The density function is then 
\begin{equation}\label{def.GPD.density}
    g_{X}(x) = \begin{cases}
               \dfrac{1}{\sigma} \Big(1 + \dfrac{\xi x}{\sigma} \Big)^{-1/\xi - 1}, &\quad  \xi \neq 0               \\\\
               
               \dfrac{1}{\sigma}e^{-x/\sigma}, &\quad  \xi = 0.\\
      
           \end{cases}
\end{equation}
The GPD is contains three distributions. When $\xi > 0 $ the GPD  is an Pareto distribution of the second kind (or Lomax distribution in the insurance literature) with a heavy tail decaying at a polynomial speed; when $\xi =0$ the GPD is an exponential distribution with a medium tail decaying exponentially; for $\xi < 0 $, the GPD becomes a short-tailed distribution the upper bound of the distribution support is finite. The $k$th moment of the GPD is existent for $\xi <1/k$; for instance, the mean and variance are finite only when $\xi <1$ and $\xi <0.5$, respectively.\\

Now we define the exponentiated GPD (exGPD) as the logarithm of the GPD random variable. That is, when $X$ is  GPD distributed, random variable $Y=\log X $ is said to be an exGPD random variable. A simple algebra gives its df as 
\begin{align}
\nonumber F_{Y}(y)&=P(Y \leq y ) = P(\log X \leq y) = P(X \leq e^{y})\\
\label{def.exGPD.df1}&=G_{X}(e^{y})=1 - \Big(1 + \dfrac{\xi e^{y}}{\sigma} \Big)^{-1/\xi} ; \quad \xi \neq 0,
\end{align}
of which the support is $-\infty < y < \infty$ for $\xi > 0$, and $-\infty < y \leq \log (-\sigma/\xi)$ for $\xi < 0$. When $\xi = 0$, we have
\begin{equation}
\label{df.exGPD.zero_xi}
F_{Y}(y)=1 - e^{-e^{y}/\sigma} = 1 - \exp (-e^{y-\log\sigma }); \quad \xi = 0, \, -\infty < y < \infty,
\end{equation} which is the Type III extreme value distribution (or the distribution of $-Y$ is that of Gumbel) with location parameter $\log\sigma$. Combining these,  we can write the density of the exGPD as
\begin{equation}\label{def.exGPD.density}
    f_{Y}(y) = \begin{cases}
               \dfrac{e^{y}}{\sigma} \Big(1 + \dfrac{\xi e^{y}}{\sigma}\Big)^{-1/\xi -1}, & \xi \neq 0               \\\\
               
               \dfrac{1}{\sigma}e^{y-e^{y}/\sigma}, & \xi = 0.\\
      
           \end{cases}
\end{equation}
We note that the role of $\sigma$ changes from a scale parameter under the GPD to a location parameter under the exGPD, as seen from $e^{y}/\sigma = e^{y-\log \sigma}$ in (\ref{def.exGPD.density}) for all $\xi$. In what follows we denote $X$ and to be the GPD variable, and $Y$ for the exGPD to avoid confusion.\\

An alternative way to create the exGPD is to define $Y$ through $Y= \log (W-d)$ conditional on $W>d$, where $W$ is the Pareto random variable with df $F_W(w)=1-(w/\beta)^{-\alpha}$, $w>\beta$. The survival function of $Y$ is  then given by
\begin{align}
\nonumber P(Y>y)   &= P( \log (W-d)>y |W>d)  =\displaystyle \frac{P(W-d>e^y|W>d)}{P(W>d)} \\
\label{}    &  =\displaystyle \frac{P(W>d+e^y)}{P(W>d)}=\frac{\big{(}\frac{d+e^y}{\beta} \big{)}^{-\alpha}}{\big{(}\frac{d}{\beta} \big{)}^{-\alpha}}= \Big{(} 1+\frac{e^y}{d}\Big{)}^{-\alpha},
\end{align} which corresponds to exGPD$(d/\alpha, 1/\alpha)$.\\

\begin{figure}[H]
\centering
\includegraphics[scale=0.20]
{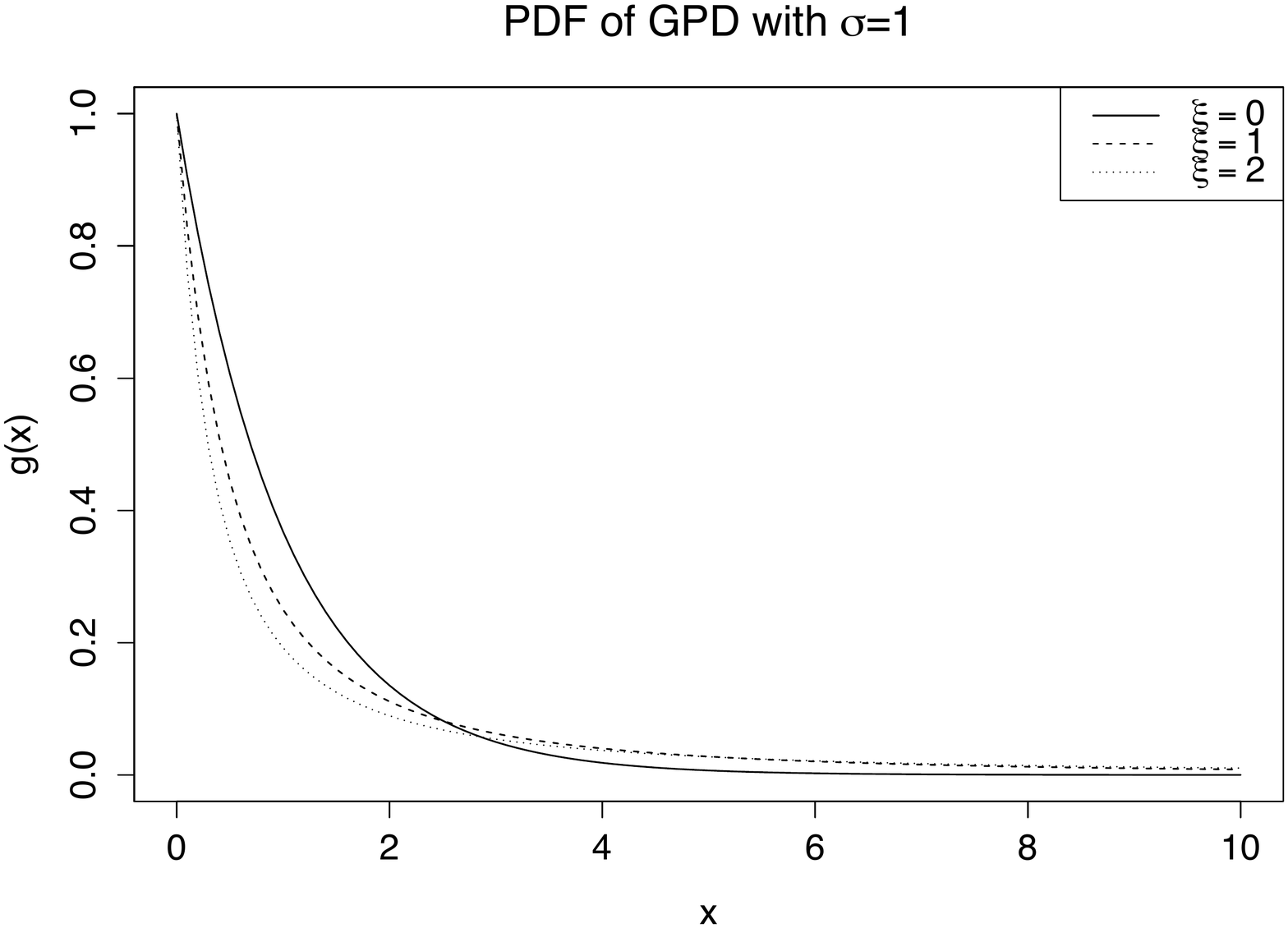}
\includegraphics[scale=0.20]{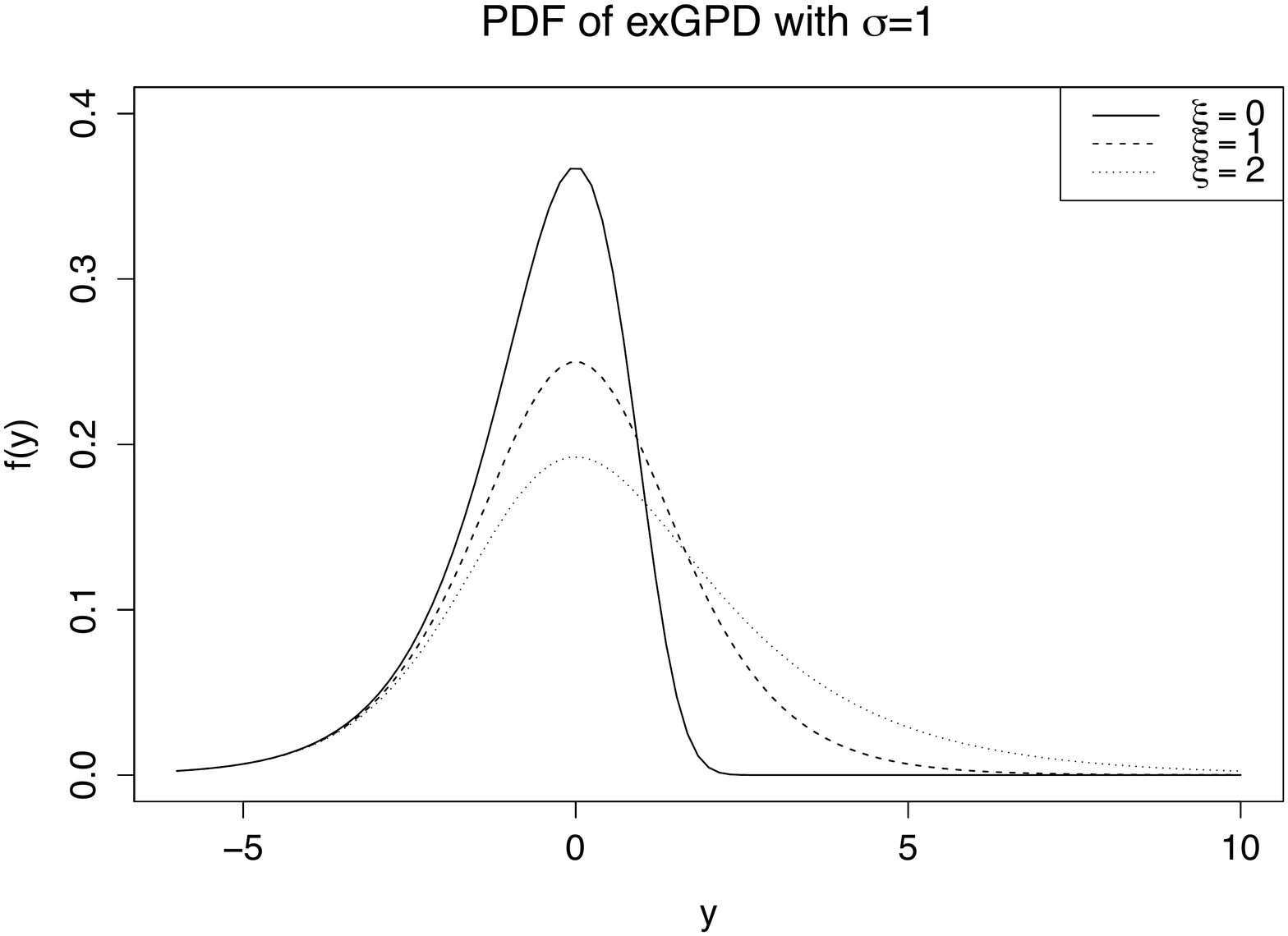}
%
\includegraphics[scale=0.20]
{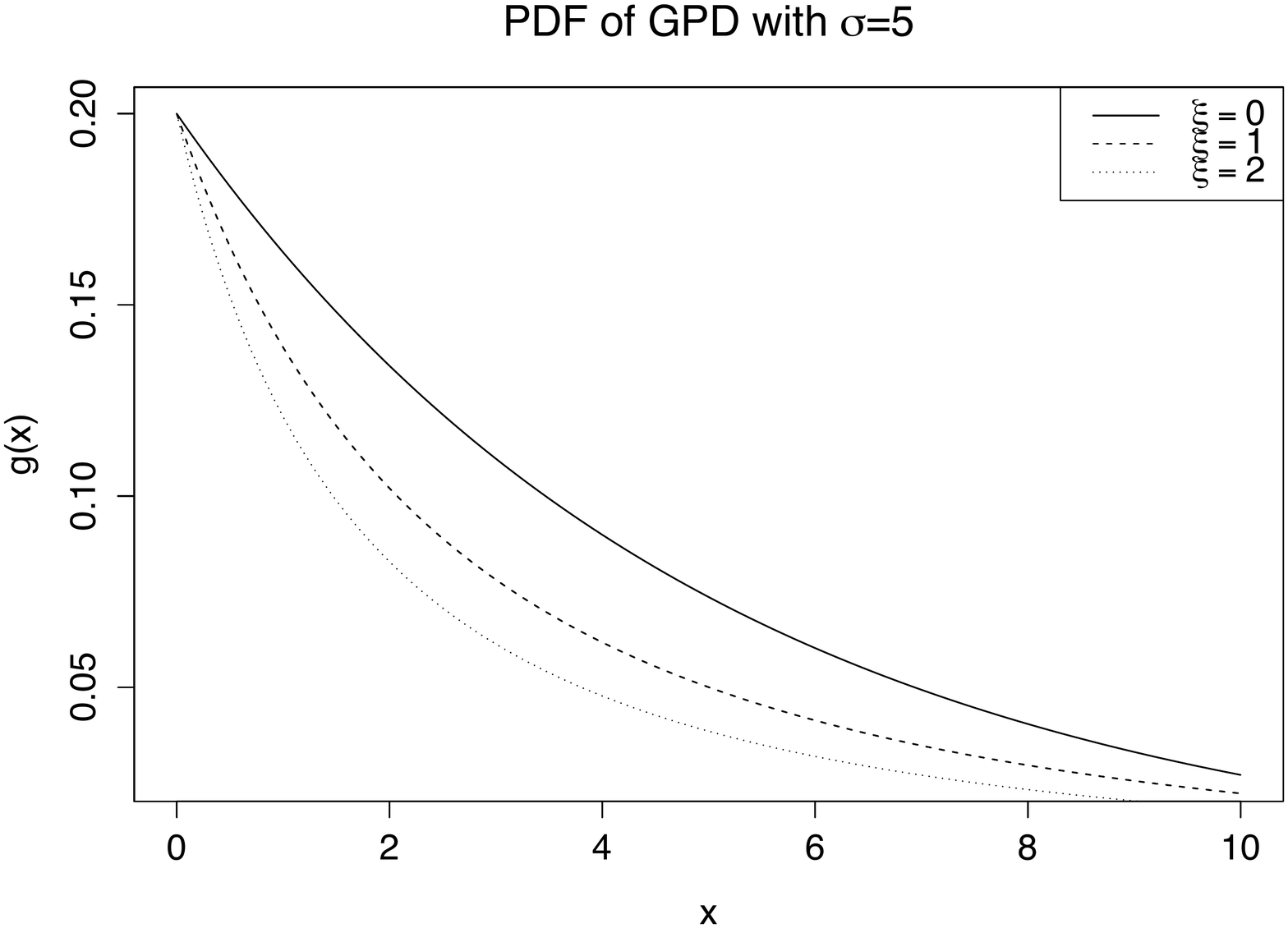}
\includegraphics[scale=0.20]{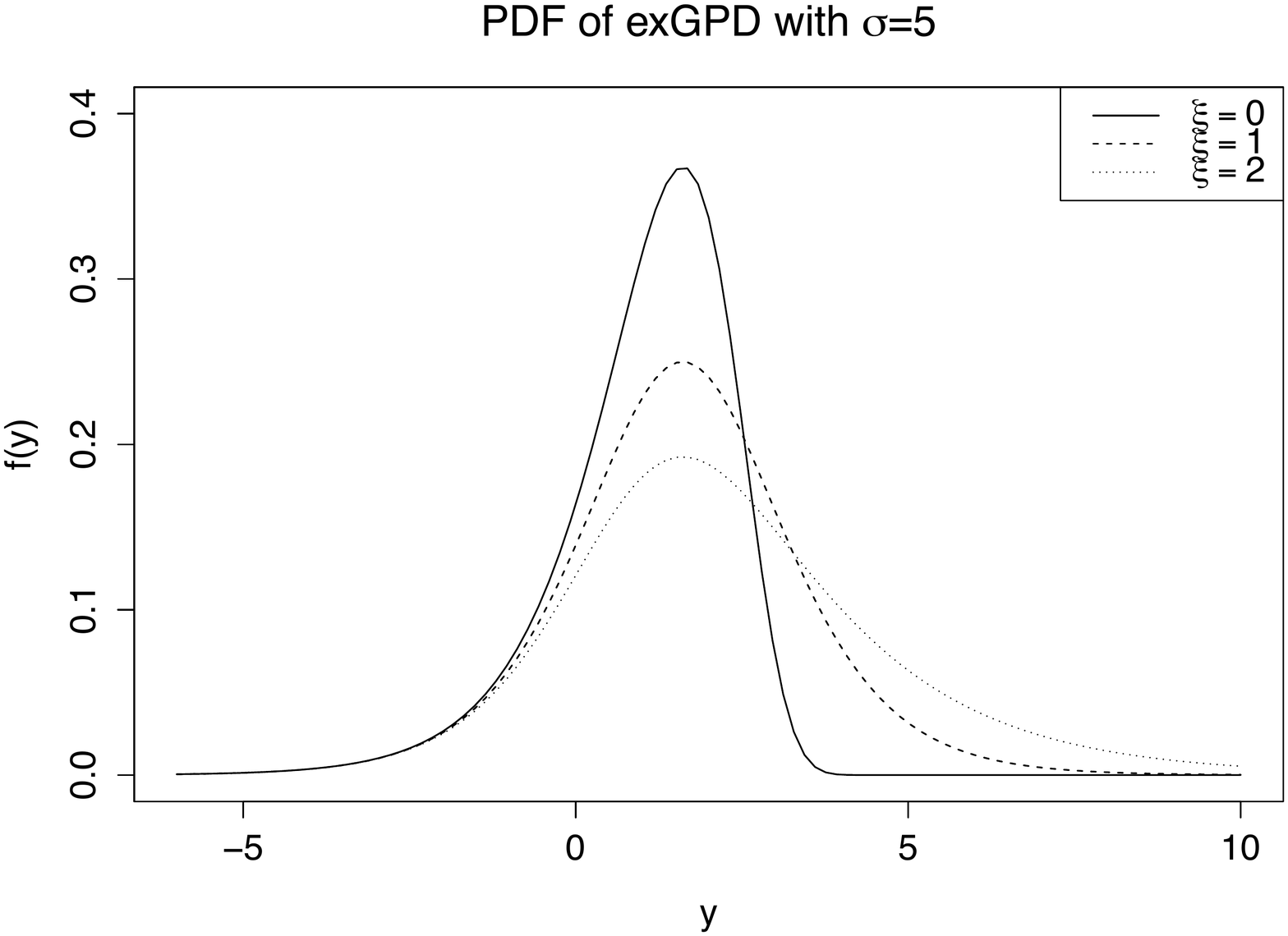}
\caption{
Density comparison: GPD vs. exGPD for $\xi \ge 0$}
\label{fig:PDFCOMPARISON}
\end{figure}

In Figure \ref{fig:PDFCOMPARISON} we compare the densities of the GPD and exGPD side by side for selected parameter choices for $\xi \ge 0$. The most dramatic change is the shape itself, where the GPD is always decreasing on its support $(0, \infty)$ whereas the exGPD, defined on the entire real line, has a peak in the middle of the distribution. Obviously, all the realized GPD values between 0 and 1 are mapped to negative numbers under the exGPD. 
Also, as seen from the densities (\ref{def.GPD.density}) and (\ref{def.exGPD.density}), the polynomial right tail of the GPD changes to an exponentially decaying tail under the exGPD through log transformation. Later we will prove the moments of all orders are actually finite for the exGPD. From the figure, we have the following additional comments:
\begin{itemize}
  \item As $\sigma$ increases the density of the exGPD shifts to the right because $\log \sigma$ is the location parameter. In fact, the mode can be shown to be $\log\sigma$, provided it exists; the proof will be given shortly. As $\sigma$ is strictly positive, the mode can take any real values with $\sigma=1$ being the boundary of the sign of the mode.
  \item As $\xi$ increases the exGPD density becomes less peaked and shows larger dispersion or variance. This implies that the shape parameter $\xi$ of the GPD roughly plays the role of a scale parameter under the exGPD. This is also hinted from the density (\ref{def.exGPD.density}); if the constant 1 is removed from the density function, $\xi$ becomes a scale parameter.
  \item The visual advantage of the exGPD over the GPD is clear from the figure. The exGPD expresses the tail thickness, represented by $\xi$, in a much clearer way in that the tail decaying angle becomes steeper as $\xi$ gets smaller. Thus one can quickly investigate how heavy the tail is by looking at, e.g., the histogram of the log data, which is less straightforward with the GPD densities.
\end{itemize}

\begin{figure}[H]
\centering
\includegraphics[scale=0.20]{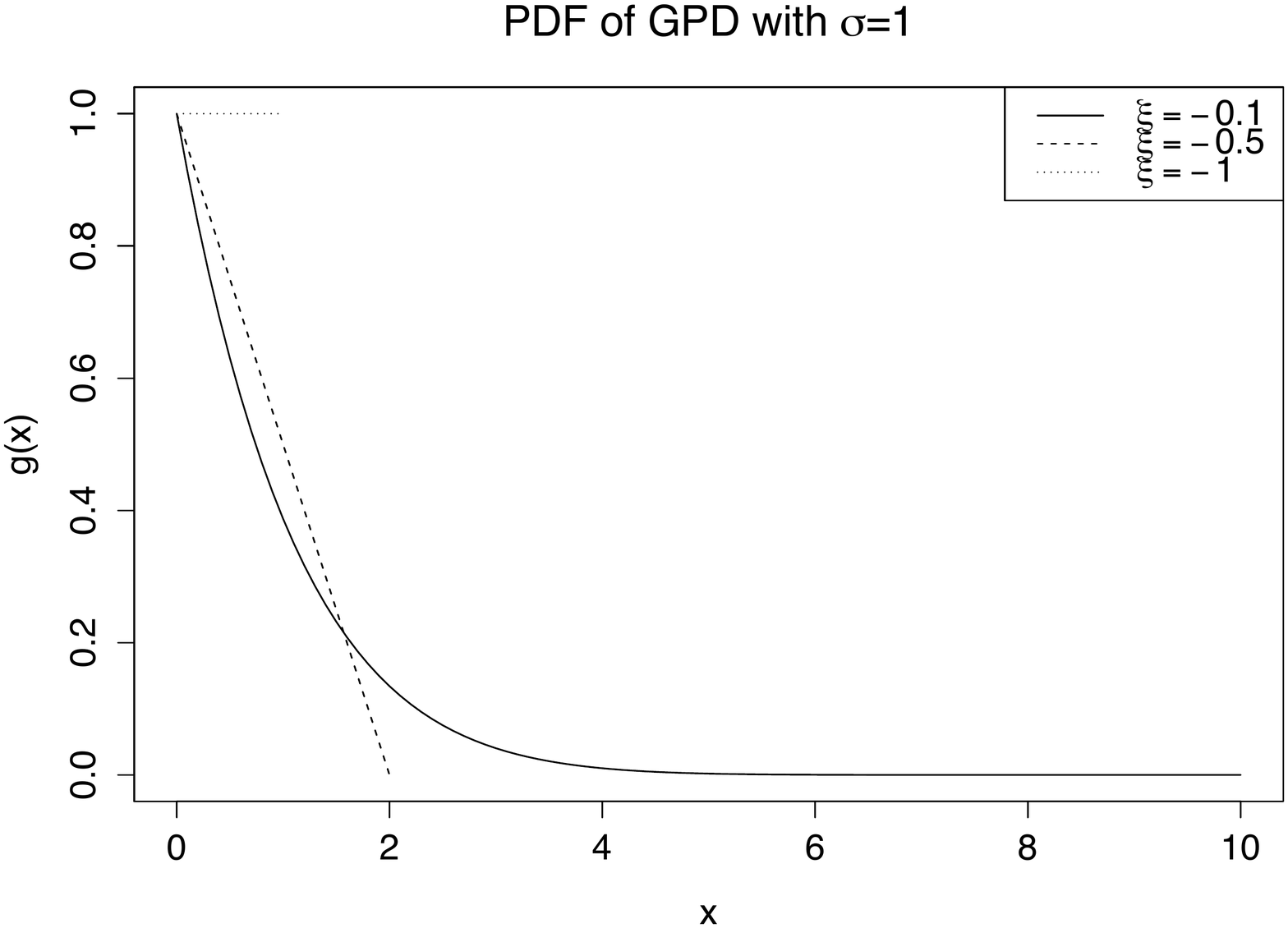}
\includegraphics[scale=0.20]{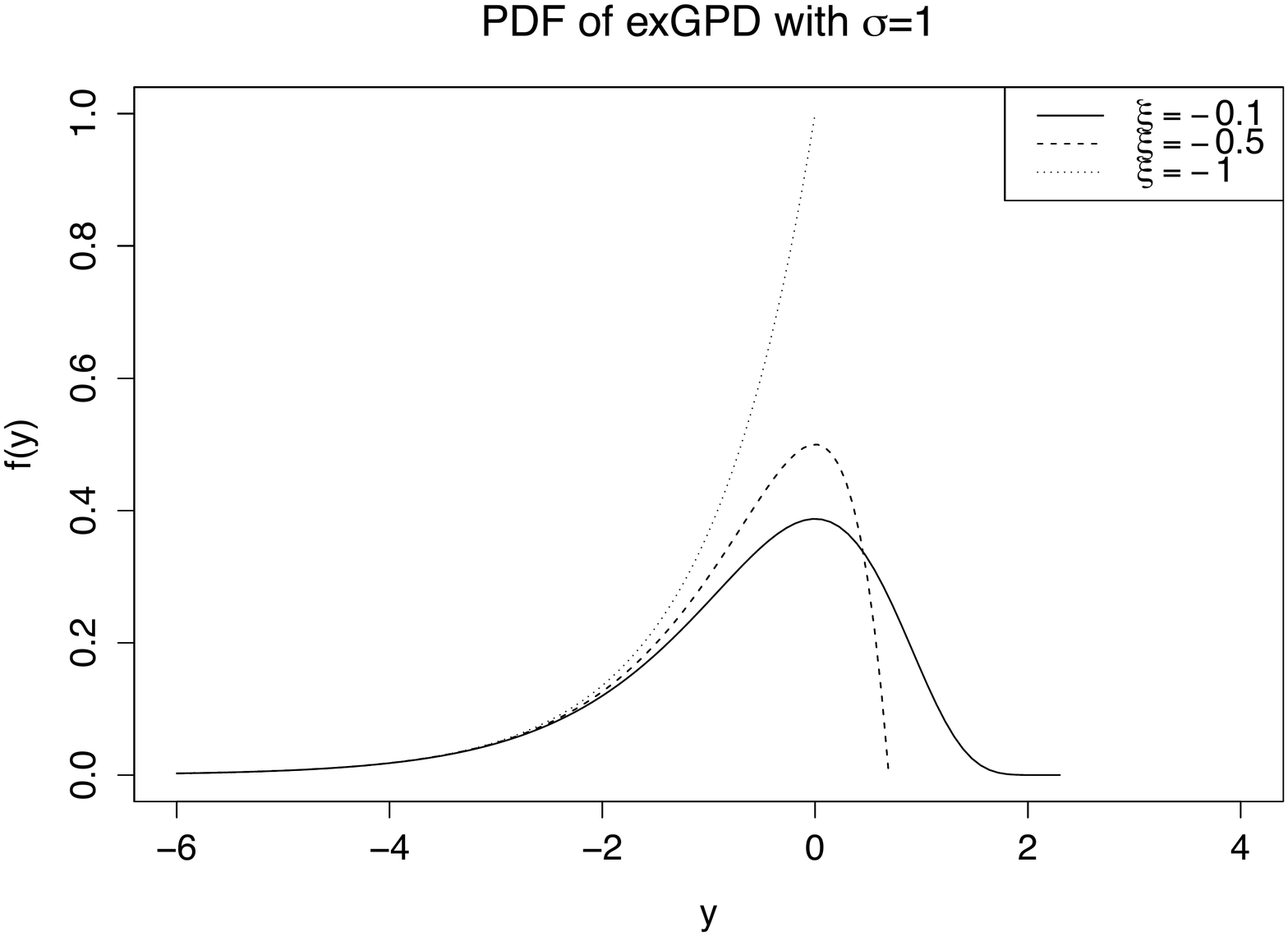}
\includegraphics[scale=0.20]{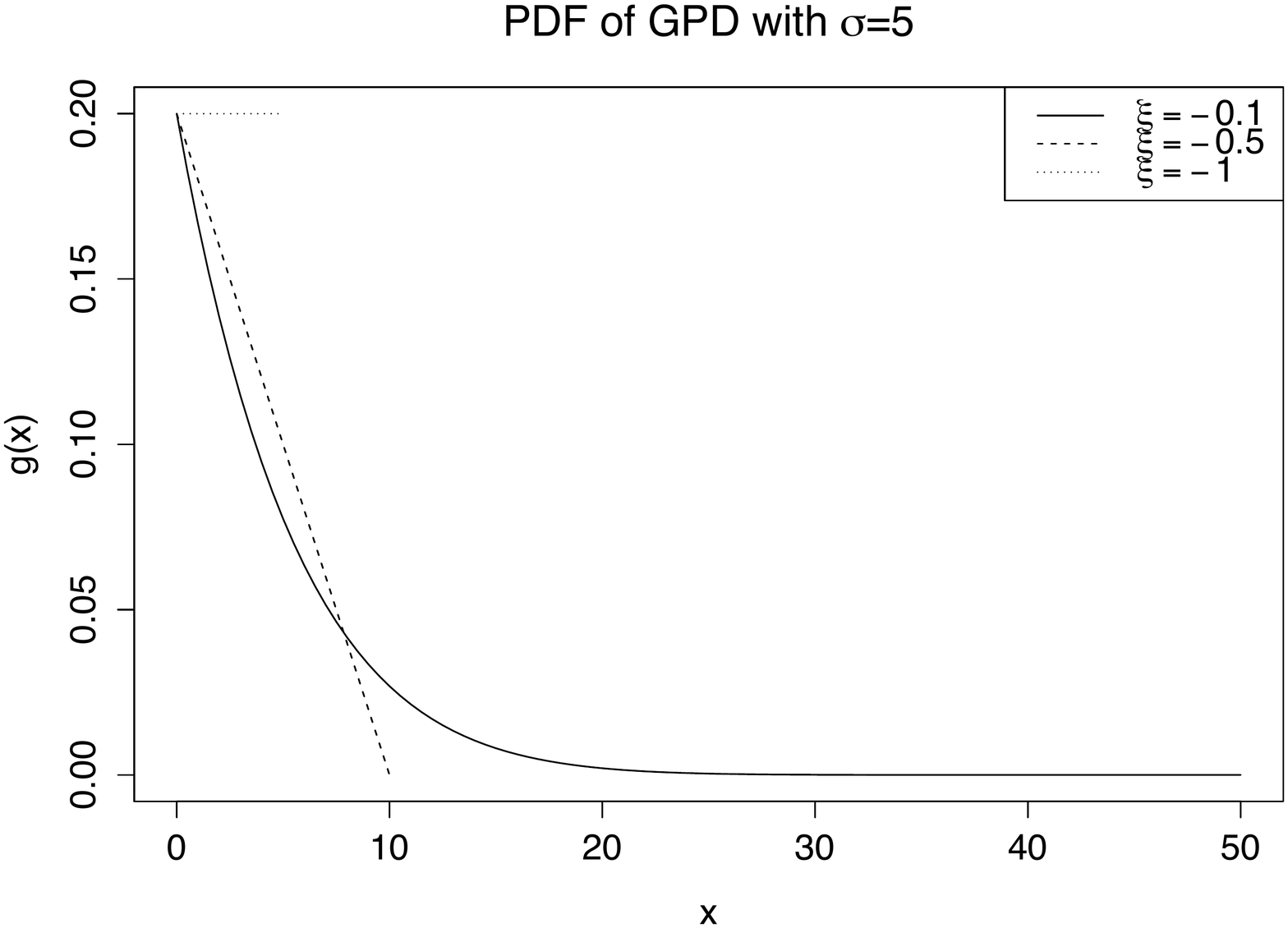}
\includegraphics[scale=0.20]{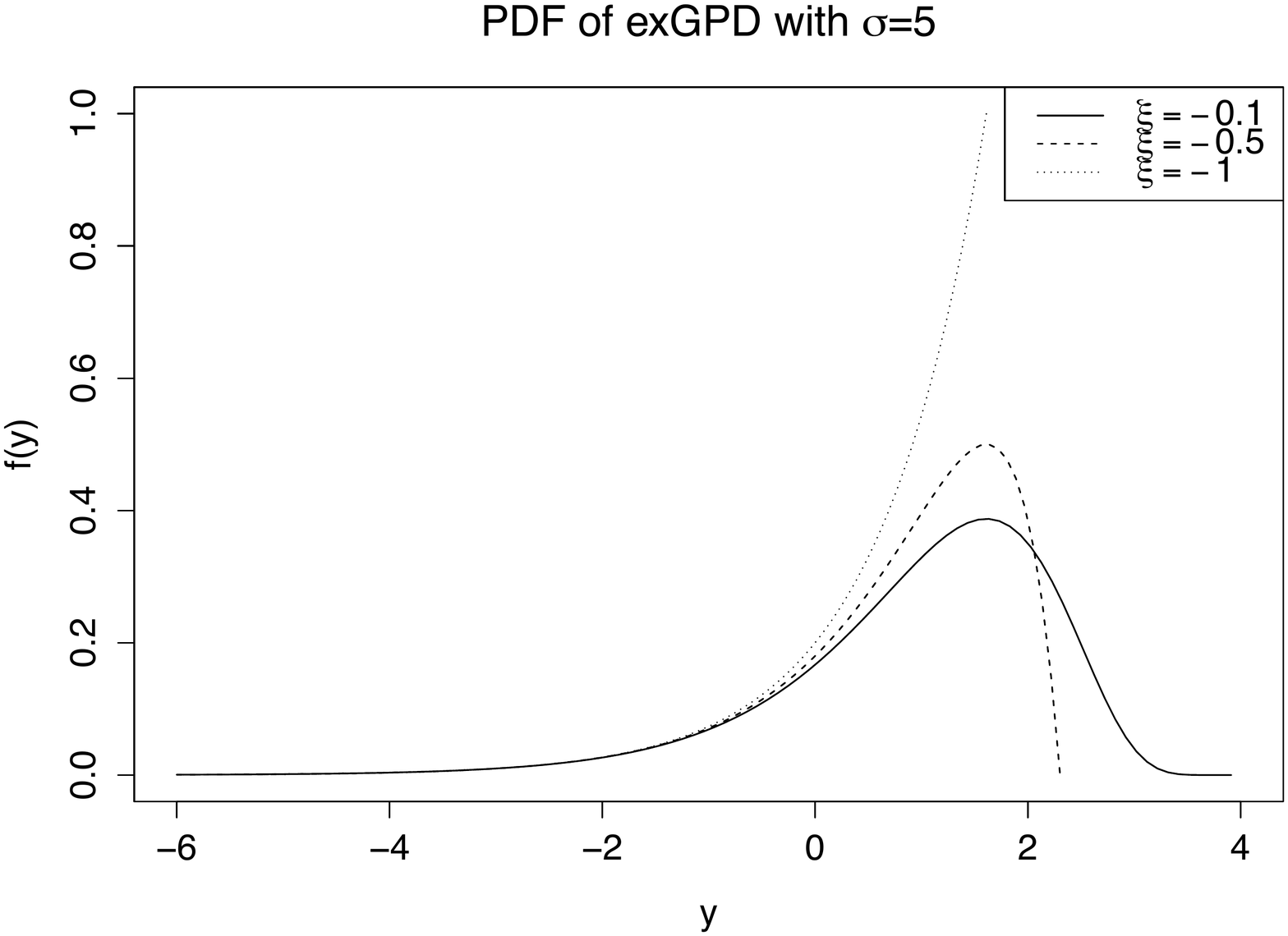}
\caption{Density comparison: GPD vs. exGPD for $\xi <0$}
\label{fig:PDFCOMPARISON.neg.xi}
\end{figure}

In Figure \ref{fig:PDFCOMPARISON.neg.xi}, the densities of the GPD and exGPD are compared for $\xi<0$. Again, the shape of the density changes substantially by log-transform. Note that there is an upper limit in the support of both distributions in this case. However, when $\xi$ is a small negative value, the shape looks not that different from $\xi \ge 0$ case in that its right tail gradually decays as if there is no upper limit to our eyes. As $|\xi|$ gets larger, but not greater than 1, there is no smooth landing around the upper limit of the support; the right tail abruptly drops to zero without hesitation. When $\xi=-1$, the exGPD density takes a drastically different shape, by shooting upside at its upper limit.  Thus, when the shape parameter is negative, its magnitude can completely change the shape of the exGPD density, as is the case for the GPD distribution. We note that the role of $\sigma$ is two-fold when $\xi<0$; it acts as the location parameter and also controls the upper limit of the distribution support.

The following result formally shows how the shape parameter of the exGPD relates to its density shape, which is also linked to the existence of the mode.
\begin{lem}\label{thm.exGPD.mode} The density of the exGPD($\sigma, \xi$) in (\ref{def.exGPD.density}) is bounded only for $\xi \ge -1$, in which case there is a unique mode at $\log \sigma$.
\end{lem}
\noindent \textbf{Proof:} To find the mode, a simple algebra yields that 
\begin{align}
\nonumber f'_{Y}(y) &= \dfrac{e^{y}}{\sigma}\bigg(1 + \dfrac{\xi e^{y}}{\sigma}\bigg)^{-1/\xi -1}
+
\dfrac{e^{y}}{\sigma}\bigg(-\dfrac{1}{\xi}-1 \bigg) \bigg(1 + \dfrac{\xi e^{y}}{\sigma}\bigg)^{-1/\xi -2}\cdot\dfrac{\xi e^{y}}{\sigma}\\
\nonumber &=
\dfrac{e^{y}}{\sigma}\bigg(1 + \dfrac{\xi e^{y}}{\sigma}\bigg)^{-1/\xi -1}
\bigg[1 + \bigg(-\dfrac{1}{\xi}-1 \bigg)
\bigg(
1+\dfrac{\xi e^{y}}{\sigma}
\bigg)^{-1}
\dfrac{\xi e^{y}}{\sigma}
\bigg]\\
\label{derivative.exGPD.pdf}&=f_{Y}(y) \frac{1-e^{y}/\sigma}{1+\xi e^{y}/\sigma}.
\end{align}Note that in the last expression, both $f_{Y}(y)$ and $1+\xi e^{y}/\sigma$ are non-negative regardless of the sign of $\xi$. We now investigate whether $f'_{Y}(y) =0$ gives a sensible solution for different ranges of $\xi$. \\
(1) $\xi<-1$ case: Using the distribution support $-\infty < y \leq \log (-\sigma/\xi)$ for negative $\xi$, we see that the numerator of the second term in (\ref{derivative.exGPD.pdf}) is bounded by
\begin{equation}
\label{ineq.sandwich.1}
1+1/\xi \le 1-e^{y}/\sigma <1.
\end{equation}When $\xi<-1$, the lower bound of this inequality becomes a strictly positive number, so $f'_{Y}(y) =0$ has no solution as seen in (\ref{derivative.exGPD.pdf}). In fact, the density value in this case gets indefinitely large as $y$ approaches its upper limit since
\begin{align}
\nonumber \lim_{y\rightarrow \log(-\sigma/\xi)^{-}}
f_{Y}(y)
\nonumber &=
\lim_{y\rightarrow \log(-\sigma/\xi)^{-}}
\dfrac{e^{y}}{\sigma}
\bigg(
1+\dfrac{\xi e^{y}}{\sigma}
\bigg)^{-1/\xi -1}\\
&=
\dfrac{-1}{\xi}
\lim_{y\rightarrow \log(-\sigma/\xi)^{-}} \bigg(
1+\dfrac{\xi e^{y}}{\sigma}
\bigg)^{-1/\xi -1}
=
+\infty.
\end{align} The last equality holds because $-1/\xi-1 < 0$ when $\xi < -1$. Thus the density is unbounded in this range.\\
(2) $\xi =-1$ case: The density reduces to $f_{Y}(y)=e^{y}/\sigma$, the exponential function shifted by $\log \sigma$. This function, defined on $-\infty < y \leq \log \sigma$, is increasing in $y$ and attains its maximum of 1, at its upper limit of the support, $y=\log \sigma$, which is the mode.\\
(3) $-1<\xi <0$ case: The lower bound of (\ref{ineq.sandwich.1}) is negative in this range, so $f'_{Y}(y) =0$ has a unique solution at $\log \sigma$ from (\ref{derivative.exGPD.pdf}). Note that the mode in this case is strictly inside the distribution support because $-\infty < \log \sigma < \log (-\sigma/\xi)$, in contrast to the $\xi =-1$ case where the maximum occurs at the upper boundary of the support.\\
(4) $\xi \ge 0$ case: In this case the distribution support is $-\infty <y<\infty$, which yields
\begin{equation}
\label{ineq.sandwich.2}
-\infty < 1-e^{y}/\sigma <1.
\end{equation}Thus $f'_{Y}(y) =0$ has a unique solution $y=\log \sigma$, which is  the mode. \quad $\Box$\\

\begin{figure}[H]
\centering
\includegraphics[scale=0.20]
{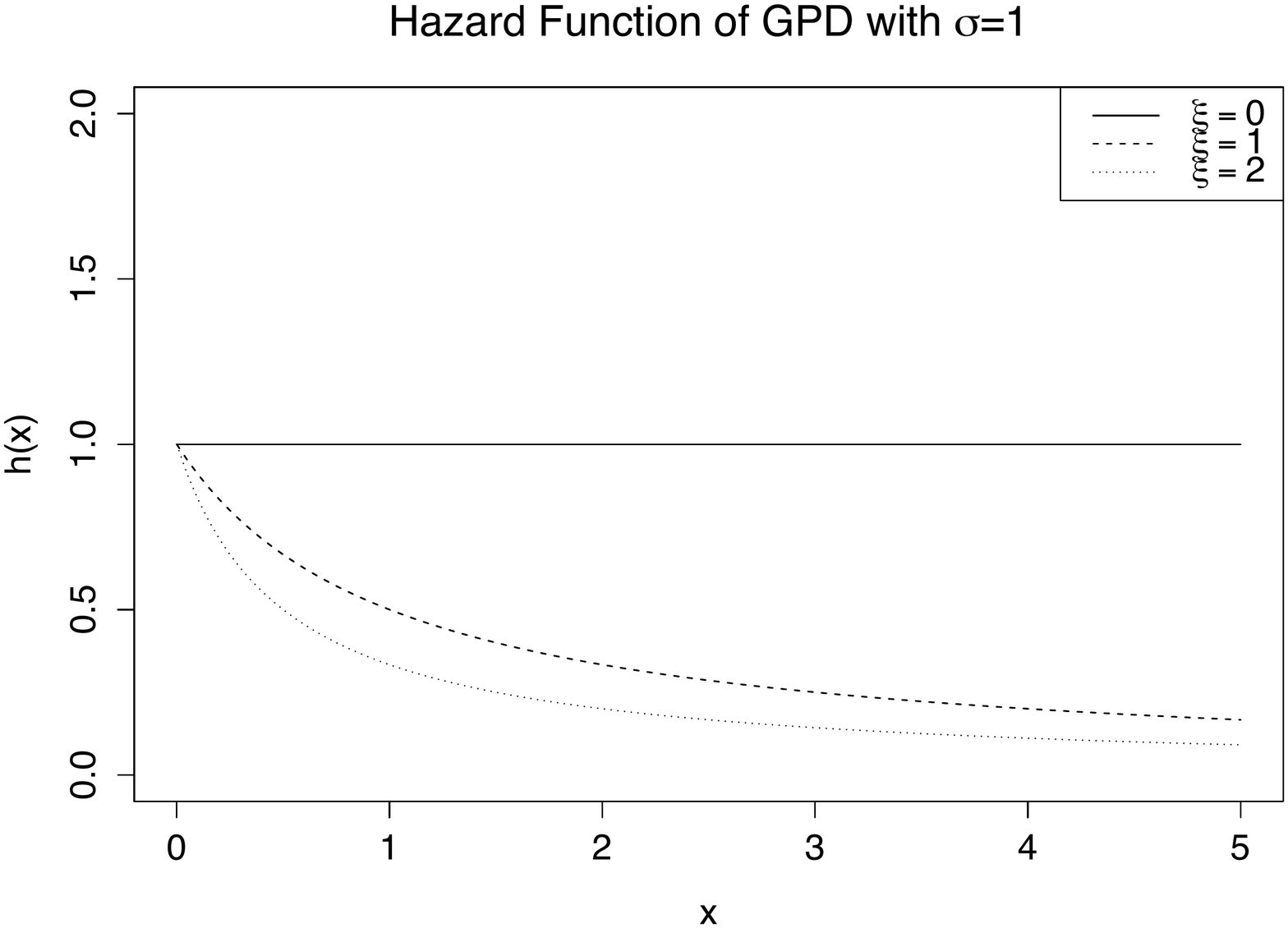}
\includegraphics[scale=0.20]{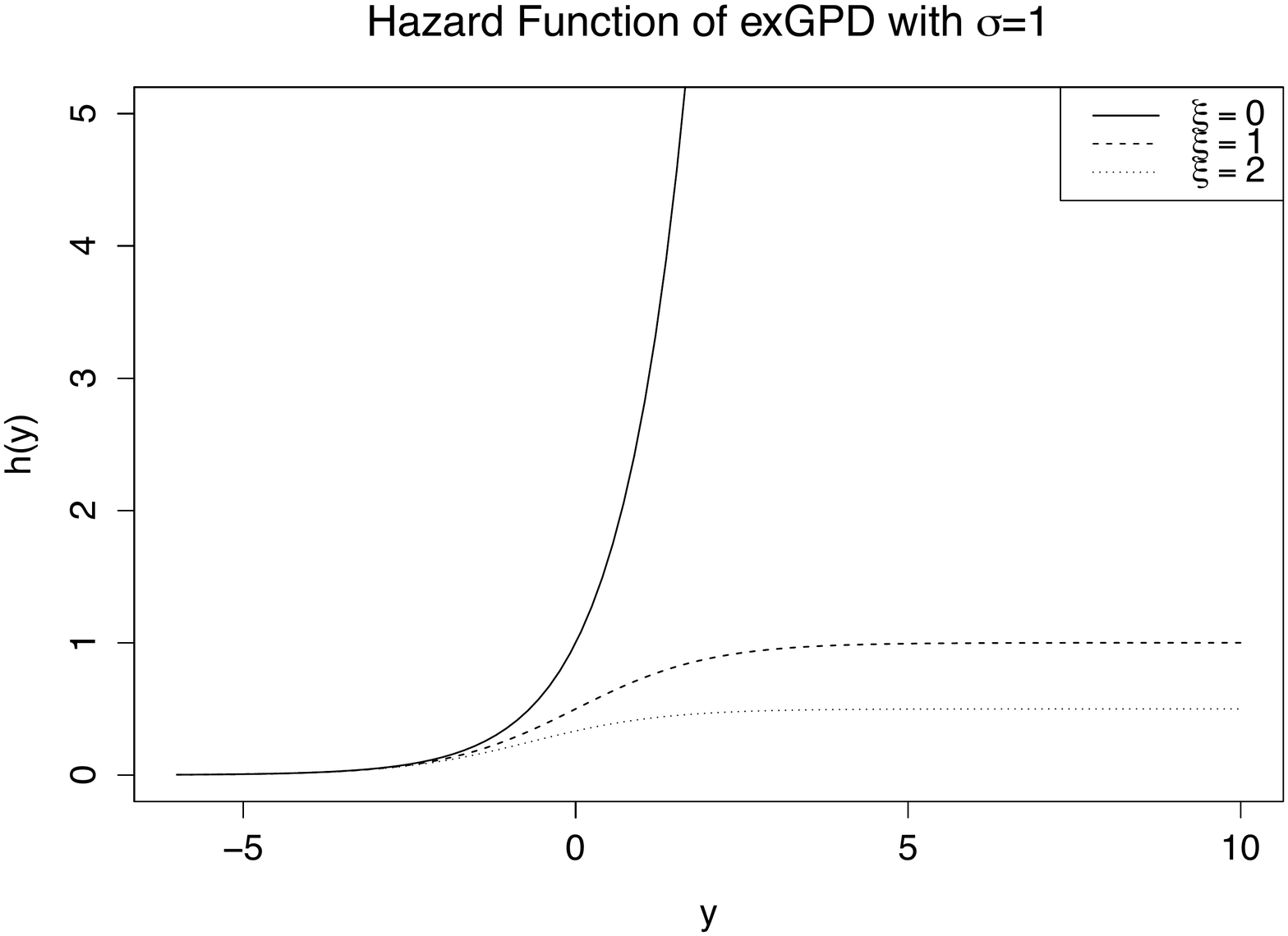}
\caption{
Hazard function comparison between GPD and exGPD with different parameters }
\label{fig:HAZARDCOMPARISON}
\end{figure}

Now we turn to the hazard function of the exGPD which is given by
\begin{equation}
\label{hazard.ft.exGPD}
h_{Y}(y)=\dfrac{f_{Y}(y)}{1-F_{Y}(y)}=
\dfrac{e^{y}}{\sigma}\bigg( 1+\dfrac{\xi}{\sigma}e^{y}\bigg)^{-1}
=\dfrac{e^{y}/\sigma}{1+\xi e^{y}/\sigma}
=\dfrac{1}{\sigma e^{-y}+ \xi}, \quad \xi \ne 0.
\end{equation}
Figure \ref{fig:HAZARDCOMPARISON} compares the hazard function of the GPD and exGPD. The hazard function of the GPD can increase or decrease depending on the sign of $\xi$, with heavy tail implied for $\xi>0$ as the hazard function decreases in $y$. However the exGPD always has an increasing hazard function, a.k.a.  increasing fairlure rate (IFR), regardless of the sign of $\xi$, as seen from (\ref{hazard.ft.exGPD}) and the figure. According to the standard theory, therefore, we see that exGPD is also DMRL (decreasing mean residual lifetime), indicating a light-tailed distribution.
\begin{figure}[H]
\centering
\includegraphics[scale=0.20]
{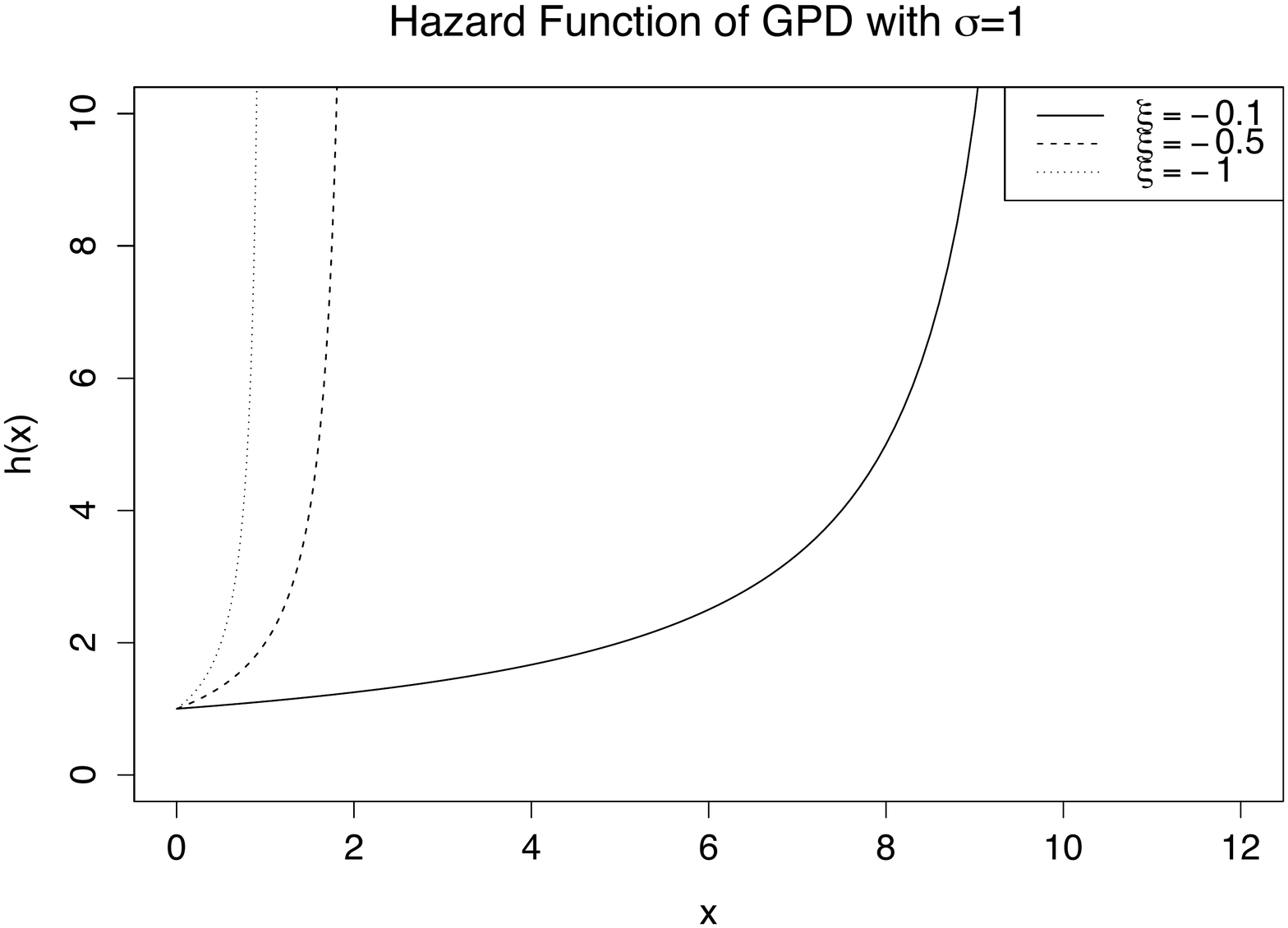}
\includegraphics[scale=0.20]{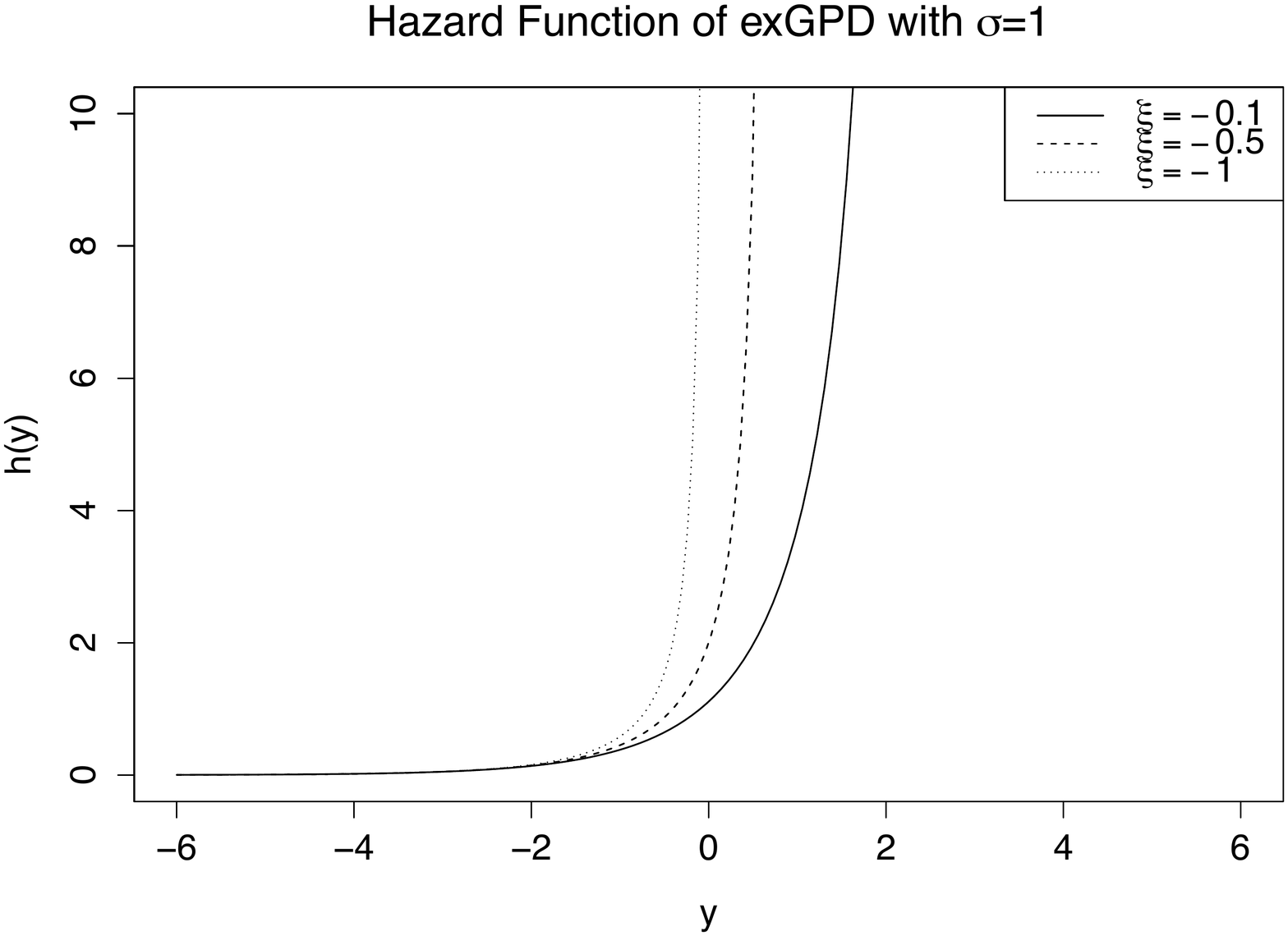}
\end{figure}

\subsection{Moment Generating Function}
In what follows, we denote $X$ to be a GPD($\sigma, \xi$) random variable and $Y=\log X $ to be the corresponding exGPD($\sigma, \xi$) variable. In order to derive the moment generating function (mgf) of $Y$, we note the following relationship:
\begin{equation}
M_{Y}(s)=E[e^{sY}]=E[e^{s \log X} ] = E[e^{\log X ^{s}} ] = E[X^{s}], \qquad s \in {\rm I\!R} . 
\end{equation}
from which we have the next result.
\begin{lem}\label{thm.exGPD.mgf} For the exGPD($\sigma, \xi$) distribution defined in (\ref{def.exGPD.df1}), the moment generating function is given by
\begin{align}\label{exGPD.mgf}
M_{Y}(s) = 
\begin{cases} 
-\dfrac{1}{\xi}
\bigg(
-\dfrac{\sigma}{\xi}
\bigg)^{s}
B(s+1,-1/\xi); \quad s \in (-1, \infty),\, \xi <0
\\\\
\dfrac{1}{\xi}
\bigg(
\dfrac{\sigma}{\xi}
\bigg)^{s}
B(s+1,1/\xi-s); \quad s \in (-1,1/\xi), \, \xi >0
\\\\
\sigma^s \Gamma (1+s); \quad \, s \in (-1, \infty), \, \xi=0, 
\end{cases}		
\end{align}
where $B$ is the beta function defined as 
\begin{equation}
B(x,y) =
\int_{0}^{1}
t^{x-1}
(1-t)^{y-1}
dt=
\dfrac{\Gamma(x)\Gamma(y)}{\Gamma(x+y)}, \quad x>0,\, y>0.
\end{equation}

 \end{lem}
\noindent \textbf{Proof:} We prove this for three cases depending on the sign of $\xi$.\\
(1) $\xi<0$ case: From the support of $X \sim $ GPD($\sigma, \xi$), we have
\begin{align}
M_{Y}(s)
&=E[X^{s}]
=\int_{0}^{-\sigma/\xi}
x^{s}\cdot
\dfrac{1}{\sigma}
\bigg(
1+\dfrac{\xi}{\sigma}x
\bigg)^{-1/\xi -1} 
dx.
\end{align}
If we let $\dfrac{\xi}{\sigma}x=-t$, we have
\begin{align}
\nonumber M_{Y}(s)
&=
\int_{0}^{1}
\bigg(
-\dfrac{\sigma}{\xi}t
\bigg)^{s}
\cdot
\dfrac{1}{\sigma}
(1-t)^{-1/\xi-1}
\cdot
-\dfrac{\sigma}{\xi}
dt
\\
\nonumber &=
-\dfrac{1}{\xi}
\bigg(
-\dfrac{\sigma}{\xi}
\bigg)^{s}
\int_{0}^{1}
t^{s}
(1-t)^{-1/\xi-1}
dt\\
&=
-\dfrac{1}{\xi}
\bigg(
-\dfrac{\sigma}{\xi}
\bigg)^{s}
B(s+1,-1/\xi); \quad -1< s.
\end{align}

\noindent (2) $\xi>0$ case:  As $X$ is defined for all positive values, 
\begin{align}
M_{Y}(s)
&=E[X^{s}]
=\int_{0}^{\infty}
x^{s}\cdot
\dfrac{1}{\sigma}
\bigg(
1+\dfrac{\xi}{\sigma}x
\bigg)^{-1/\xi -1} 
dx.
\end{align}
By letting $\dfrac{\xi}{\sigma}x=t$, we have
\begin{align}
\nonumber M_{Y}(s)
&=
\int_{0}^{\infty}
\bigg(
\dfrac{\sigma}{\xi}t
\bigg)^{s}
\cdot
\dfrac{1}{\sigma}
(1+t)^{-1/\xi-1}
\cdot
\dfrac{\sigma}{\xi}
dt
\\
\nonumber &=
\dfrac{1}{\xi}
\bigg(
\dfrac{\sigma}{\xi}
\bigg)^{s}
\int_{0}^{\infty}
t^{s}
(1+t)^{-1/\xi-1}
dt\\
&=
\dfrac{1}{\xi}
\bigg(
\dfrac{\sigma}{\xi}
\bigg)^{s}
B(s+1,1/\xi-s); \quad -1 <s <1/\xi,
\end{align}
where the last equality comes from a property of the beta function
\begin{equation}
B(x,y) =
\int_{0}^{1}
t^{x-1}
(1-t)^{y-1}
dt=
\int_{0}^{\infty}
t^{x-1}
(1+t)^{-x-y}
dt, \quad x>0, \, y>0.
\end{equation}
It is pointed out that if we compare the two mgf's ($\xi<0$ vs. $\xi>0$) in (\ref{exGPD.mgf}), only the second argument of the beta function is different. 

\noindent (3) $\xi=0$ case: The derivation is straightforward and omitted; alternatively, this can be adapted and derived from the mgf of Gumbel distribution as shown in, e.g., Ch. 22 of \cite{Johns+Kotz+etal:94b}.
  \quad $\Box$\\

\subsection{Moments}
From the mgf of the exGPD given in Lemma \ref{thm.exGPD.mgf}, we may determine the moments by differentiating it with respect to $s$. We present the first two moments for $\xi \ne 0$ first. Let us rewrite mgf in (\ref{exGPD.mgf}) as a function of gamma functions for easier differentiation.
\begin{equation}
\label{mgf.exGPD.in.gamma}
M_{Y}(s) = 
\begin{cases} 
-\dfrac{1}{\xi}
\Gamma(-1/\xi)
\cdot
\bigg(
-\dfrac{\sigma}{\xi}
\bigg)^{s}
\cdot
\dfrac{\Gamma(s+1)}{\Gamma(s-1/\xi+1)}
; \quad s \in (-1, \infty),\, \xi <0.
\\\\
\dfrac{1}{\xi}
\Gamma(1/\xi+1)
\cdot
\bigg(
\dfrac{\sigma}{\xi}
\bigg)^{s}
\cdot
\Gamma(s+1)
\cdot
\Gamma(1/\xi-s); \quad s \in (-1,1/\xi), \, \xi >0.
\end{cases}		
\end{equation}
For $\xi<0$ case, we may obtain the first two derivatives of the mgf as follows.
\begin{align*}
\dfrac{d}{ds} M_{Y}(s)
&=
-\dfrac{1}{\xi}
\Gamma(-1/\xi)
\cdot
\bigg(
-\dfrac{\sigma}{\xi}
\bigg)^{s}
\cdot
\Gamma(s+1)
\cdot
\dfrac{1}{\Gamma(s-1/\xi+1)}\\
&\quad
\quad
\cdot
\bigg(
\log \bigg( -\dfrac{\sigma}{\xi}\bigg)
+
\psi(s+1)
-
\psi(s-1/\xi+1)
\bigg).\\
\dfrac{d^{2}}{ds^{2}} M_{Y}(s)
&=
-\dfrac{1}{\xi}
\Gamma(-1/\xi)
\cdot
\bigg(
-\dfrac{\sigma}{\xi}
\bigg)^{s}
\cdot
\Gamma(s+1)
\cdot
\dfrac{1}{\Gamma(s-1/\xi+1)}\\
&
\quad
\quad
\cdot\bigg[
\bigg(
\log \bigg( -\dfrac{\sigma}{\xi}\bigg)
+
\psi(s+1)
-
\psi(s-1/\xi+1)
\bigg)^{2}
+
0+\psi'(s+1)-\psi'(s-1/\xi+1)
\bigg],
\end{align*} where $\psi$ is the digamma function. 
In this derivation, we used $\Gamma'(s+1) = \Gamma(s+1)\cdot \psi(s+1)$ and 
\[\bigg(\dfrac{1}{\Gamma(s-1/\xi+1)}\bigg)' = -\dfrac{\psi(s-1/\xi+1)}{\Gamma(s-1/\xi+1)}.\]
Similarly, when $\xi>0$, we use the fact $\Gamma'(1/\xi-s) = -\Gamma(1/\xi-s)\psi(1/\xi-s)$ to get 
\begin{align*}
\dfrac{d}{ds} M_{Y}(s)
&=
\dfrac{1}{\xi}
\Gamma(1/\xi+1)
\cdot
\bigg(
\dfrac{\sigma}{\xi}
\bigg)^{s}
\cdot
\Gamma(s+1)
\cdot
\Gamma(1/\xi-s)
\bigg(
\log \bigg( \dfrac{\sigma}{\xi} \bigg)
+
\psi(s+1)
-
\psi(1/\xi-s)
\bigg).\\
\dfrac{d^{2}}{ds^{2}} M_{Y}(s)
&=
\dfrac{1}{\xi}
\Gamma(1/\xi+1)
\cdot
\bigg(
\dfrac{\sigma}{\xi}
\bigg)^{s}
\cdot
\Gamma(s+1)
\cdot
\Gamma(1/\xi-s)\\
& 
\quad
\quad
\cdot
\bigg[
\bigg(
\log \bigg( \dfrac{\sigma}{\xi} \bigg)
+
\psi(s+1)
-
\psi(1/\xi-s)
\bigg)^{2}
+0
+\psi'(s+1)
+\psi'(1/\xi-s)
\bigg].
\end{align*}
Hence, by setting $s=0$, we have the first moment of the exGPD:
\begin{align}
\label{mean.exGPD}E[Y] = 
\dfrac{d}{ds} M_{Y}(s)
\bigg|_{s=0}
=
\begin{cases} 
\log
\bigg(
-\dfrac{\sigma}{\xi}
\bigg)
+\psi(1)
-
\psi
(1-1/\xi)
;\quad \xi <0.
\\\\
\log
\bigg(
\dfrac{\sigma}{\xi}
\bigg)
+\psi(1)
-
\psi
(1/\xi)
;\quad \xi >0.\\ \\
\log \sigma +\psi(1); \quad \xi=0.
\end{cases}		
\end{align}
Here $\psi(1)=-\gamma$, where $\gamma=0.5772 \cdots$ is the Euler-Mascheroni constant. The case for for $\xi=0$ has been determined separately, but it is easily done. 
Likewise, the second moment is given by
\begin{align}
\nonumber E[Y^{2}] = 
\dfrac{d^{2}}{ds^{2}} M_{Y}(s)
\bigg|_{s=0}
&=
\begin{cases} 
\bigg(
\log
\bigg(
-\dfrac{\sigma}{\xi}
\bigg)
+\psi(1)
-
\psi
(1-1/\xi)
\bigg)^{2}
+\psi'(1)
-
\psi'(-1/\xi+1)
;\quad \xi <0
\\\\
\bigg(
\log
\bigg(
\dfrac{\sigma}{\xi}
\bigg)
+\psi(1)
-
\psi
(1/\xi)
\bigg)^{2}
+\psi'(1)
+
\psi'(1/\xi)
;\quad \xi >0
\\\\
\big(\log \sigma + \psi (1)\big)^{2} + \psi'(1);
\quad \xi =0.
\end{cases}
\\ 
\nonumber\\
&=
\begin{cases} 
E[Y]^{2}
+
\psi'(1)
-
\psi'(-1/\xi+1)
;\quad \xi <0
\\\\
E[Y]^{2}
+\psi'(1)
+
\psi'(1/\xi)
;\quad \xi >0
\\\\
E[Y]^{2} + \psi'(1);
\quad \xi =0.
\end{cases}		
\end{align}
From the second raw moment, the variance of the exGPD becomes
\begin{equation}
\label{var.exGPD.trigamma}Var[Y]=
\begin{cases} 
\psi'(1)
-
\psi'(-1/\xi+1)
\\\\
\psi'(1)
+
\psi'(1/\xi)
\\\\
\psi'(1)
\end{cases}
=\begin{cases} 
\dfrac{\pi^{2}}{6}
-
\displaystyle \sum_{k=1}^{\infty}
\dfrac{1}{(k-1/\xi)^{2}}
;\quad \xi <0
\\
\dfrac{\pi^{2}}{6}
+
\displaystyle \sum_{k=1}^{\infty}
\dfrac{1}{(k+1/\xi-1)^{2}}
;\quad \xi >0
\\
\dfrac{\pi^{2}}{6};
\quad \xi =0.
\end{cases}		
\end{equation}
Note that the variance of exGPD depends only on $\xi$. Furthermore, since the summation terms in (\ref{var.exGPD.trigamma}) are always positive, the value of the sample variance of $Y$, denoted $s_Y^2$, can serve as a quick indicator about the sign of $\xi$. That is, if $s_Y^2>\pi^2/6 \simeq 1.645$, then $\xi$ is deemed positive; otherwise, $\xi$ is negative. The first two moments yield the method of moments estimator (MME) of the exGPD parameter. For $\xi<0$,
\begin{align*}
\hat{\xi}_{MME}	
&=
\dfrac{1}{1-(\psi')^{-1}(\psi'(1)-s^{2})},\\
\hat{\sigma}_{MME}	
&=
-\hat{\xi}_{MME}	
\cdot
e^{\bar{Y} -\psi(1)+\psi(1-1/\hat{\xi}_{MME}	)}
\end{align*}
and for $\xi >0$,
\begin{align*}
\hat{\xi}_{MME}	
&=
\dfrac{1}{(\psi')^{-1}(s^{2}-\psi'(1))},\\
\hat{\sigma}_{MME}
&=
\hat{\xi}_{MME}	
\cdot
e^{\bar{Y} -\psi(1)+\psi(1/\hat{\xi}_{MME}	)}.
\end{align*}
When $\xi =0$, we have $\hat{\sigma}_{MME} = e^{\bar{Y} -\psi(1)}. $\\

%
One may further obtain higher moments by differentiating the mgf repeatedly, but the task is not straightforward as higher order derivatives can be complicated, involving infinite series, as seen from the derivation of the first two moments above. 
However, formally establishing the existence of higher moments of the exGPD is important because the difficulty of handling the GPD, such as sampling variability, is essentially attributed to its tail heaviness, directly connected to the non-existence of its moments. 
\begin{cor}
The exGPD defined in (\ref{def.exGPD.density}) has finite moments of all orders.
\end{cor} 
\noindent \textbf{Proof:} 
Referring to (\ref{mgf.exGPD.in.gamma}), the mgf of the exGPD is a product of three functions in terms of $s$, for both $\xi>0$ and $\xi<0$ cases. The first term is simply an exponential function of $s$, which has derivatives of all orders. The second and third terms are gamma functions or its reciprocal. The derivatives of  a gamma function can be written through the polygamma function. The  polygamma function of order $k$ is  defined as the $(k+1)$-th derivative of the logarithm of the gamma function
\begin{equation}
\label{ }
\psi^{{(k)}}(z)=\frac{d^{k}}{dz^{k}}\psi(z)=\frac{d^{k+1}}{dz^{k+1}}\log \Gamma(z)=(-1)^{k+1}k! \sum^{\infty}_{{r=0}}\frac{1}{(z+r)^{k+1}}, \quad k \ge 1,
\end{equation}
which is finite for $z>0$. Hence we conclude that the mgf of the exGPD has derivatives of all orders, each of which in turn gives a finite value when evaluated at $0$. A similar argument holds for $\xi=0$ case.
 \qquad $\Box$

\subsection{Further properties}
In the GPD literature, various distributional quantities are available in addition to the ordinary moments and MEF; see,e.g., Ch 3 of \cite{EmbrechtsModellingExtremalEvents}. Here we list similar properties of the exGPD; some are parallel to those of the GPD, but others are different. We present the findings first, and then provide the proofs. All findings include both $\xi>0$ and $\xi<0$ cases.\\

\noindent (a) For a real value $r$ with $r\xi > -1$, 
\begin{equation}
\label{ }
E\bigg[\bigg(1+\dfrac{\xi}{\sigma}e^{Y}\bigg)^{-r}\bigg] 
= \dfrac{1}{1+r\xi} ; \quad r\xi > -1
\end{equation}
\noindent (b) For a non-negative integer $k$,
\begin{equation}
\label{ }
E\bigg[\bigg(\log\bigg(1+\dfrac{\xi}{\sigma}e^{Y}\bigg)\bigg)^{k}\bigg] 
=\xi^{k}
\cdot
k!
\end{equation}
\noindent (c) For real value $r$ with $1+r>|\xi|$, 
\begin{equation}
\label{ }
E\bigg[e^{Y}
\cdot
\big(
\overline{F}_{Y}(Y)
\big)^{r}\bigg] 
= \dfrac{\sigma}{(r+1-\xi)(r+1)}
\end{equation}
\noindent 
(d) Assume that $N \sim Poi(\lambda)$, independent of the iid sequence $(Y_{i})_{i=1}^{n}$ where $Y_{i} \sim exGPD(\sigma,\xi)$, $i =1,2,\cdots n$. Write $M_{N}=\max(Y_{1},Y_{2},\cdots,Y_{n})$. Then

\begin{equation}
\label{ }
P(M_{N} \leq y)
=\exp
\bigg(
-\lambda
\bigg(
1+\dfrac{\xi e^{y}}{\sigma}
\bigg)^{-1/\xi}
\bigg)
\end{equation}
\noindent (e) For a constant $c$, 
\begin{equation}
\label{ }
\int_{c}^{\infty}
\overline{F}_{Y}(y)dy
=\int_{c}^{\infty}\bigg(1+\dfrac{\xi}{\sigma}e^{y}\bigg)^{-1/\xi}dy= \begin{cases}
B\Big((1+ \xi e^{c}/\sigma)^{{-1}}; \,1/\xi,0\Big),      & \xi>0, \\ \\
B\Big(1 + \xi e^{c}/\sigma;\,1-1/\xi,0\Big),      & \xi<0,
\end{cases}
\end{equation} where $B(x;a,b)$ is the incomplete beta function
\begin{equation}
\label{ }
B(x;a,b)
=\int_{0}^{x}
t^{a-1}
(1-t)^{b-1}
dt.
\end{equation}

\noindent \textbf{Proof:} 
Proof for (a): 
\begin{align}
\nonumber E\bigg[\bigg(1+\dfrac{\xi}{\sigma}e^{Y}\bigg)^{-r}\bigg] 
&= \int_{-\infty}^{\infty}
\bigg(1+\dfrac{\xi}{\sigma}e^{y}\bigg)^{-r}
\cdot
\dfrac{e^{y}}{\sigma}
\cdot
\bigg(1+\dfrac{\xi}{\sigma}e^{y}\bigg)^{-1/\xi-1}
dy\\
\label{pro.(a).exGPD1} &=
 \int_{-\infty}^{\infty}
\dfrac{e^{y}}{\sigma}
\cdot
\bigg(1+\dfrac{\xi}{\sigma}e^{y}\bigg)^{-1/\xi-r-1}
dy.
\end{align}
For $\xi>0$ case, we take $1+\xi e^{y}/\sigma=t$, then by  integration by substitution, (\ref{pro.(a).exGPD1}) is equal to

\begin{equation}
\label{ }
\int_{1}^{\infty}
t^{-1/\xi - r -1}
\cdot
\dfrac{1}{\xi}
dt
=
\dfrac{1}{\xi}
\cdot
\dfrac{1}{1/\xi + r}
=
\dfrac{1}{1+r\xi} ; \quad r > -1/\xi.
\end{equation}
When $\xi<0$, (\ref{pro.(a).exGPD1}) becomes, with the same substitution,
\begin{equation}
\label{ }
\int_{1}^{0}
t^{-1/\xi - r -1}
\cdot
\dfrac{1}{\xi}
dt
=
\dfrac{1}{\xi}
\cdot
\dfrac{1}{1/\xi + r}
=
\dfrac{1}{1+r\xi} ; \quad r < -1/\xi,
\end{equation}
completing the proof. Note that the different conditions for $\xi$ has been combined to $ r\xi > -1$. We comment thus that the condition used in Theorem 4.3.13 of \cite{EmbrechtsModellingExtremalEvents} is not correct.\\

\noindent Proof for (b):
\begin{align}
\label{pro.(a).exGPD2} 
E\bigg[\bigg(\log\bigg(1+\dfrac{\xi}{\sigma}e^{Y}\bigg)\bigg)^{k}\bigg] 
&= \int_{-\infty}^{\infty}
\bigg(\log\bigg(1+\dfrac{\xi}{\sigma}e^{y}\bigg)\bigg)^{k}
\cdot
\dfrac{e^{y}}{\sigma}
\cdot
\bigg(1+\dfrac{\xi}{\sigma}e^{y}\bigg)^{-1/\xi-1}
dy.
\end{align}
For $\xi>0$ case, if we take $\log(1+\xi e^{y}/\sigma)=t$, then by integration by substitution, (\ref{pro.(a).exGPD2}) is equal to
\begin{equation}
\label{ }
\xi^{{-1}}\int_{0}^{\infty}
t^{k}
\cdot
e^{-t/\xi}
dt = \xi^{k}\cdot k!.
\end{equation}
For $\xi<0$ case, we obtain the same result with the same substitution, but via a slightly different integration as $\xi<0$ and $t<0$.
Note that depending on $k$ being odd or even, the quantity can be negative or positive for $\xi<0$.\\

\noindent Proof for (c): 
\begin{align}
\nonumber E\bigg[e^{Y}
\cdot
\big(
\overline{F}_{Y}(Y)
\big)^{r}\bigg] 
&=
\int_{-\infty}^{\infty}
e^{y}
\cdot
\bigg(
\bigg(1+\dfrac{\xi}{\sigma}e^{y}\bigg)^{-1/\xi}
\bigg)^{r}
\cdot
\dfrac{e^{y}}{\sigma}
\cdot
\bigg(1+\dfrac{\xi}{\sigma}e^{y}\bigg)^{-1/\xi-1}
dy\\
\label{pro.(a).exGPD3} &=
\int_{-\infty}^{\infty}
e^{y}
\cdot
\dfrac{1}{\sigma}
\bigg(1+\dfrac{\xi}{\sigma}e^{y}\bigg)
^{-(r+1)/\xi -1}
dy.
\end{align}
For $\xi>0$ case, if we take $1+\xi e^{y}/\sigma=t$, then (\ref{pro.(a).exGPD3}) is equal to 
\begin{align}
\nonumber \int_{1}^{\infty}
\dfrac{\sigma}{\xi}
\cdot
(t-1)
t^{-(r+1)/\xi -1}
\cdot
\dfrac{1}{\xi}
dt
\nonumber &=
\dfrac{\sigma}{\xi^{2}}
\int_{1}^{\infty}
(
t^{-(r+1)/\xi}
-
t^{-(r+1)/\xi -1}
)
dt\\
\nonumber &=
\dfrac{\sigma}{\xi^{2}}
\bigg(
\dfrac{\xi}{r+1-\xi}
-
\dfrac{\xi}{r+1}
\bigg)\\
&=
\dfrac{\sigma}{(r+1-\xi)(r+1)};\quad 1+r>\xi.
\end{align}
For $\xi<0$ case, with the same substitution, (\ref{pro.(a).exGPD3}) becomes
\begin{equation}
\label{ }
\int_{1}^{0}
\dfrac{\sigma}{\xi}
\cdot
(t-1)
t^{-(r+1)/\xi -1}
\cdot
\dfrac{1}{\xi}
dt
=\dfrac{\sigma}{(r+1-\xi)(r+1)};\quad r+1>0.
\end{equation} Thus the conditions differ depending on the sign of $\xi$ even though the results are identical. If one wishes to merge these for convenience sake, then $1+r>|\xi|$ would serve the purpose.\\

\noindent  Proof for (d): Note that 
\begin{align*}
P(M_{N} \leq y)
=E [P( M_{N}\leq y | N=n )]
=
\sum_{n=0}^{\infty}
P(M_{N}\leq y|N=n)\cdot P(N=n).
\end{align*}
On the other hand, we have
\begin{align*}
P(M_{N}\leq y|N=n)
&=P(M_{n}\leq y)=P(Y_{1}\leq y,Y_{2}\leq y,\cdots,Y_{n}\leq y)\\
&= P(Y_{1}\leq y) P(Y_{2}\leq y)\cdots P(Y_{n}\leq y)=(F_{Y}(y) )^{n}.
\end{align*} 
Therefore, we have
\begin{align}
\nonumber P(M_{N} \leq y)
&=
\sum_{n=0}^{\infty}
(F_{Y}(y) )^{n}\cdot \dfrac{\lambda^{n} e^{-\lambda}}{n!}
=
e^{-\lambda}
\sum_{n=0}^{\infty}
\dfrac{(\lambda \cdot F_{Y}(y) )^{n}}{n!}\\
\nonumber &=
e^{-\lambda}
\cdot
e^{\lambda \cdot F_{Y}(y)}
= e^{-\lambda(1-F_{Y}(y))}
= e^{-\lambda \cdot \overline{F}_{Y}(y)}\\
&=\exp
\bigg(
-\lambda \,
\bigg(
1+\dfrac{\xi e^{y}}{\sigma}
\bigg)^{-1/\xi}
\bigg).
\end{align}
This property is related to the POT framework. If transformed back, this is the generalized extreme value (GEV) distribution.\\

\noindent  Proof for (e): For $\xi>0$ case, use integration by substitution $t=(1+\xi e^{y}/\sigma)^{-1}$ to get
\begin{align}
\nonumber \int_{c}^{\infty}
\overline{F}_{Y}(y)dy
&=
\int_{c}^{\infty}
\bigg(
1+\dfrac{\xi e^{y}}{\sigma}
\bigg)^{-1/\xi}
dy
=\int_{0}^{(1+\xi e^{c}/\sigma)^{-1}}
t
^{1/\xi-1}
\cdot
(1-t)^{-1}
dt\\
&=
B\Big((1 + \xi e^{x}/\sigma)^{-1};1/\xi,0\Big),
\end{align} which can be evaluated using the hypergeometric function.
Similarly, for $\xi<0$ case, we set $t=1+\xi e^{y}/\sigma$ to get 
\begin{align}
\nonumber \int_{c}^{\log(-\sigma/\xi)}
\overline{F}_{Y}(y)dy
&=
\int_{c}^{\log(-\sigma/\xi)}
\bigg(
1+\dfrac{\xi e^{y}}{\sigma}
\bigg)^{-1/\xi}
dy=\int_{1+\xi e^{c}/\sigma}^{0}
t
^{-1/\xi}
\cdot
(t-1)^{-1}
dt \\
&=\int^{1+\xi e^{c}/\sigma}_{0}
t
^{-1/\xi}
\cdot
(1-t)^{-1}
dt=
B\Big(1 + \xi e^{c}/\sigma;\, 1-1/\xi,0\Big).
\end{align}
For $\xi =0$ case, we may use integration by substitution $t=e^{y-\log\sigma}$ to get 
\begin{equation}
\label{int.Fbar.zero_xi}
\int_{c}^{\infty}
\overline{F}_{Y}(y)dy
=
\int_{c}^{\infty}
e^{-e^{y-\log \sigma}}
dy=
\int_{e^{c-\log\sigma}}^{\infty}
t^{-1} e^{-t}
dt=\Gamma(0,e^{c-\log \sigma}),
\end{equation}
 where $\Gamma(s,x)$ is the incomplete gamma function
\begin{equation}
\label{ }
\Gamma(s,x) 
=\int_{x}^{\infty}
t^{s-1}
e^{-t}
dt. \quad \quad \Box
\end{equation}
Sometimes the three-parameter GPD is found in the literature, which is created by adding a location parameter $\mu$ to the GPD. The df of this distribution is defined as $G(x; \mu, \sigma, \xi)=G(x-\mu; \sigma, \xi)$, $x>\mu$, where $G(x;\sigma, \xi)$ is the GPD in (\ref{def.GPD.df1}). We comment that most results in this section can be readily applied to the three-parameter GPD without additional difficulty.

\subsection{Maximum likelihood estimation}
Let $y_1, \cdots, y_n$ be an iid sample from exGPD$({\sigma}, {\xi})$. Then from its density (\ref{def.exGPD.density}) the log-likelihood of the exGPD can be written as
\begin{equation}
l(\sigma, \xi) = \sum_{i=1}^{n} y_{i} + \dfrac{n \log \sigma}{\xi} -\Big(\dfrac{1}{\xi}+1 \Big)\sum_{i=1}^{n} \log (\sigma + \xi e^{y_{i}} ).
\end{equation}
and the MLE $(\hat{\sigma}, \hat{\xi})$ may be found from solving 
\begin{align*}
\dfrac{\partial l}{\partial \sigma}&=\dfrac{n}{\sigma \xi} 
- \Big(\dfrac{1}{\xi} +1\Big) \sum_{i=1}^{n}\Big(\dfrac{1}{\sigma+\xi e^{y_{i}}} \Big) 
=0      \\
\dfrac{\partial l}{\partial \xi}
&=-\dfrac{n \log \sigma}{\xi^{2}}
- \Big(\dfrac{1}{\xi} +1\Big) \sum_{i=1}^{n}\Big(\dfrac{e^{y_{i}}}{\sigma+\xi e^{y_{i}}} \Big)
+ \dfrac{1}{\xi^{2}}\sum_{i=1}^{n}\log (\sigma + \xi e^{y_{i}})
=0.
\end{align*} 
Clearly, the MLE has to be solved numerically. Due to the condition that $\sigma+\xi e^{y_{i}}>0$ for all $i$, one must have $y_{(n)} < \log (-\sigma/\xi)$ for  $\xi < 0$ where $y_{(n)}$ is the sample maximum.  Then, for any given $\xi < -1$, we see that
\begin{align}
\lim_{\log(-\sigma/\xi) \rightarrow y_{(n)}^{+}} l(\sigma, \xi) = \infty,
\end{align} 
implying that there is no finite MLE. So the MLE exists only when $\xi >-1$,  the  same condition for the MLE of the GPD to exist as pointed out in, e.g., \cite{grimshaw1993computing}. In fact, it turns out that the MLEs for the GPD and exGPD are identical for a given sample. To see this, consider the densities of these two distributions with the same parameter $(\sigma, \xi)$, that is, $g_X(x)$ and $f_Y(y)$  from (\ref{def.GPD.density}) and (\ref{def.exGPD.density}), respectively. When a sample $x_1, \cdots, x_n$ is given from the GPD, we can obtain a corresponding exGPD sample $y_1, \cdots, y_n$ with $y_i=\log x_i$. The log-likelihood function for the exGPD ($\xi \ne 0$) is then shown from the density to be
\begin{align}
\nonumber \sum_{i=1}^n \log f_Y(y_i|\sigma, \xi)   &=  \sum_{i=1}^n y_i +\sum_{i=1}^n\log g_X(e^{y_i}|\sigma, \xi) \\
\label{}    &  = \sum_{i=1}^n \log x_i +\sum_{i=1}^n\log g_X(x_i|\sigma, \xi) 
\end{align}
In the last expression, the first term involves no parameter and the second term stands for the log-likelihood function of the GDP. Hence maximizing the likelihood functions for both distributions leads to the identical parameter estimate. This implies that, even though the log transform may stablize the volatility of the GPD sample in terms of, e.g., moments, the MLE does not offer additional benefits from this stabilization. Consequently, if we let 
$\hat{\theta} 
=
\begin{pmatrix}
\hat{\sigma}\\
\hat{\xi}
\end{pmatrix}$ be the MLE of the parameter 
$\theta
=
\begin{pmatrix}
\sigma\\
\xi
\end{pmatrix}
$ for exGPD$(\sigma, \xi)$, we have
\begin{equation}
\sqrt{n}(\hat{\theta} - \theta) \overset{d}{\sim} BVN( \mathbf{ 0},[I_{1}(\theta)]^{-1}),
\end{equation}
where the covariance matrix is given by 
\begin{equation}
\label{ }
[I_{1}(\theta)]^{-1}=
\begin{bmatrix}
2 \sigma^{2}(1+\xi)&
-\sigma(1+\xi)\\
-\sigma(1+\xi)& 
(1+\xi)^{2}\\
\end{bmatrix}, \quad \xi >-0.5.
\end{equation} 
The condition $\xi >-0.5$ is needed for the information matrix to be defined properly; see, e.g., \cite{Smith:84a} or \cite{EmbrechtsModellingExtremalEvents}.


\section{Risk measures and mean excess function}
Modern financial risk management in the EVT framework often requires determination of some well-known tail risk measures to summarize the riskiness of the underlying loss distribution. Two popular such measures are the Value-at-Risk (VaR) and the conditional tail expectation (CTE). The VaR of a continuous random variable $Y$ at the 100$p\%$ level, with df $F_Y$, is the 100$p$ quantile of the distribution of $Y$, denoted by
\begin{equation}
y_p=F_Y^{-1}(p). \label{VaR}
\end{equation}
The VaR is a widely-accepted standard risk measure used in solvency and risk analyses in financial and insurance industry, with $p$ close to 1. From (\ref{def.exGPD.df1}), the VaR of the exGPD is easily shown to be 
\begin{equation}\label{def.exGPD.quantile}
y_{p} = \log\bigg(\dfrac{\sigma}{\xi} \Big( (1-p)^{-\xi}-1\Big)\bigg), \quad 0 < p < 1,
\end{equation} or, using that quantiles are preserved under a monotone transform, we can also write $y_p=\log x_p$, where $x_p$ is the 100$p$ quantile of the GPD.
Despite of its popularity, the quantile VaR measure has also been criticized for failing to meet one of the properties that any desirable risk measure is expected to satisfy, the criteria known as the coherent risk measure axioms; see \cite{Artz+Delb+etal:99a}. In this regard, the CTE\footnote{a.k.a. Tail VaR, Conditional VaR and Expected Shortfall in the literature.} has received much attention as an alternative,  coherent tail risk measure. The idea of the CTE is to measure the average severity of the loss when the extreme loss does occurs, where the extreme loss is represented by the VaR. The CTE of the random variable $Y$ at $100p\%$ level is defined as
\begin{equation}\label{def.cte}
CTE_p(Y)=E(Y|Y>y_p).
\end{equation}
Though the CTE is relatively new, it is closely related to the excess variable $Y-u|Y>u$ which has long been used in statistics and reliability studies to examine the tail behaviour of a given distribution. Its expectation $e(u)=E(Y-u|Y>u)$, known as the mean excess function (MEF) or the mean residual lifetime, is connected to the CTE through
\begin{equation}
\label{ }
CTE_p(Y)=E(Y|Y>y_p)=y_p+E(Y-y_p|Y>y_p)=y_p+ e(y_p).
\end{equation}
Thus the CTE is immediate if we can obtain the MEF. \\

For the exGPD we first derive the df of the excess variable, denoted $Y'=Y-u|Y>u$ for a constant $u$. Referring to the df of the exGPD (\ref{def.exGPD.df1}), we have 
\begin{align}
\nonumber F_{Y'}(y')
&=
P(Y' < y')
= P(Y-u \leq y' | Y > u) = \dfrac{P(u < Y < u + y')}
{P(Y>u)}
= \dfrac{\overline{F}_{Y}(u) - \overline{F}_{Y}(u+y')}
{\overline{F}_{Y}(u)}\\
\label{df.excess.variable.exGPD} 
&=
1-\dfrac{\bigg( 1+ \dfrac{\xi e^{u+y'}}{\sigma} \bigg)^{-1/\xi}}{\bigg( 1+ \dfrac{\xi e^{u}}{\sigma} \bigg)^{-1/\xi}}
= 1- \bigg( \dfrac{ \sigma + \xi e^{u+y'}}{ \sigma + \xi e^{u}} \bigg)^{-1/\xi}, \quad y'\ge 0, \,\,\, \xi \ne 0.
\end{align}
Unfortunately this is not an exGPD df. So, unlike the GPD case, the exGPD does not have the stability property. Nonetheless the MEF can be obtained as the mean of the df (\ref{df.excess.variable.exGPD}). We note that the expectation of any nonnegative continuous random variable $X$ can be written as
\begin{align}
E[X]
=
\int_{0}^{\infty}
x f_{X}(x)dx
=
-x (1-F_{X}(x))\Big{|}_{0}^{\infty}
+
\int_{0}^{\infty}
1-F_{X}(x)dx
=
\int_{0}^{\infty}
\overline{F}_{X}(x)dx.
\end{align} The last equality holds as long as $\lim_{x \rightarrow \infty}x\overline{F}_{X}(x)=0$. For the excess variable $Y'$ for the exGPD, which is nonnegative, this condition is met for both $\xi<0$ and $\xi>0$; the former case has a finite upper limit, and for the latter case the tail $\overline{F}_{Y'}(y')$ decays exponentially as seen from (\ref{df.excess.variable.exGPD}). Therefore the MEF of the exGPD becomes
\begin{align}
\nonumber e_Y(u)
&=
E[Y']
=E[Y-u|Y>u]
=
\int_{0}^{\infty}
\bigg(
\dfrac{\sigma + \xi e^{u} e^{y'}}{\sigma + \xi e^{u}}
\bigg)^{-1/\xi}
dy'\\
\nonumber &=
\bigg(
\dfrac{\sigma + \xi e^{u}}{\sigma}
\bigg)^{1/\xi}
\cdot
\int_{0}^{\infty}
\bigg(
\dfrac{\sigma + \xi e^{u} e^{y'}}{\sigma}
\bigg)^{-1/\xi}
dy'\\
\label{mef.exGPD.1}  &=
\bigg(
\dfrac{\sigma + \xi e^{u}}{\sigma}
\bigg)^{1/\xi}
\cdot
\int_{0}^{\infty}
\bigg(
1+
\dfrac{\xi e^{y'}}{\sigma e^{-u} }
\bigg)^{-1/\xi}
dy'
\end{align}
To evaluate the integration term in the last expression, let us consider an exGPD variable $Y^*$ with parameter $(\sigma^{*},\xi)$ where $\sigma^{*}=\sigma e^{-u} >0$. Then 
We recognize the the integrand in (\ref{mef.exGPD.1}) as the survival function of another exGPD variable $Y^*$, so that using Property (e) of the previous section,
\begin{align}
\nonumber 
   e_Y(u)&= \bigg(
\dfrac{\sigma + \xi e^{u}}{\sigma}
\bigg)^{1/\xi} \, \int_{0}^{\infty} \overline{F}_{Y^*}(y^*) dy^*\\
\nonumber &=\begin{cases}
\bigg(
1+\dfrac{ \xi e^{u}}{\sigma}
\bigg)^{1/\xi} B\Big((1+ \xi /\sigma^{*})^{{-1}}; \,1/\xi,0\Big),      & \xi>0, \\ \\
\bigg(
1+\dfrac{ \xi e^{u}}{\sigma}
\bigg)^{1/\xi} B\Big(1 + \xi /\sigma^{*};\,1-1/\xi,0\Big),      & \xi<0,
\end{cases}\\
&=\begin{cases}
(\bar{F}_{Y}(u))^{-1} B\Big((1+ \xi e^{u}/\sigma)^{{-1}}; \,1/\xi,0\Big),      & \xi>0, \\ \\
(\bar{F}_{Y}(u))^{-1} B\Big(1 + \xi e^{u}/\sigma;\,1-1/\xi,0\Big),      & \xi<0,
\end{cases}
\end{align} 
For $\xi=0$ case, the distribution of $Y'=Y-u|Y>u$ is similarly obtained as
\begin{equation}
\label{ }
F_{Y'}(y') =1-\dfrac{\overline{F}_{Y}(u+y')}{\overline{F}_{Y}(u)}
=1-e^{-e^{u+y'-\log \sigma} + e^{u-\log \sigma}}
=
1-
e^{-e^{u-\log \sigma}(e^{y'}-1)}, \quad y'\ge0,
\end{equation} and its MEF is given by
\begin{align}
\nonumber 
e(u)    &=E[Y']   =
\int_{0}^{\infty}
e^{-e^{u+y'-\log \sigma} + e^{u-\log \sigma}}dy'=
\exp (e^{u-\log \sigma})\int_{0}^{\infty}
e^{-e^{u+y'-\log \sigma}}dy'\\
\label{mef.exGPD.1} &=\exp (e^{u-\log \sigma})\int_{u}^{\infty}
e^{-e^{y'-\log \sigma}}dy'=\exp (e^{u-\log \sigma})\, \Gamma(0,e^{u}/\sigma),
\end{align} where the last equality comes from (\ref{int.Fbar.zero_xi}).

\section{An application: Finding tail index}
In this section we illustrate how the proposed exGPD distribution can be used in the EVT framework. In particular, we focus on finding the tail thickness, characterized as the tail index $\alpha$, for a given dataset. Here $\alpha$ corresponds to $1/\xi$ in the GPD. Determining the tail index is a problem of great importance as it characterizes the degree of tail risk  by assigning the given distribution to a proper maximum domain of attraction, as well as its moment existence range.  To formally motivate, let us denote the tail (or survival) function of a distribution by $\bar{F}(x)=1-F(x)$, $-\infty < x < \infty$. Then we say that $\bar{F}(x)$ is regularly varying with index $-\alpha <0$, and write $\bar{F} \in \mathcal{R}_{-\alpha}$, if
\beq\label{regular.varying.ft} \lim_{x\rightarrow \infty} \frac{\bar{F}(x\lambda)}{\bar{F}(x)}= \lambda^{-\alpha}, \qquad \lambda>0. \eeq
When $\alpha=0$ the tail is called slowly varying, or $\bar{F} \in \mathcal{R}_{0}$. Using this we can represent a regularly varying distribution as
$\bar{F}(x) \sim L(x) \, x^{-\alpha}$  where ${L}(\cdot)\in \mathcal{R}_{0}$. We note that $f(x) \sim g(x)$ means $\lim_{x \to \infty} f(x)/g(x) =1$.  Thus the tail of regularly varying functions can be represented by power functions multiplied by slowly varying functions. In many applications it is of interest to accurately estimate the tail index $\alpha$ from the given dataset; we refer the reader to, e.g., \cite{EmbrechtsModellingExtremalEvents} and \cite{Beirl+Goege+etal:06a} for further details. \\

The method of \cite{Hill:75a} is one of the most studied methods to estimate the tail index $\alpha$. When $F(x)$ is a regularly varying distribution with $\alpha >0$, the Hill method is designed to find $\alpha $ using a sample $X_1,...,X_n$ from $F(x)$. If we let the order statistics be $X_{n,n} \le ... \le X_{1,n}$, the Hill estimator is defined as
\begin{equation}
\label{hill}
\hat{\alpha}_{k,n} = \bigg{(} \frac{1}{k} \sum^k_{j=1} \log X_{j,n}-\log X_{k,n} \bigg{)}^{-1}, \quad 2 \le k \le n.
\end{equation}
After computing $\hat{\alpha}_{k,n}$ for different $k$ values, one can draw the Hill plot $\{(k, \hat{\alpha}^{-1}_{k,n})\}$, $k=2,3,...,n$ on the plane. The goal is then to find an area where $\hat{\alpha}^{-1}_{k,n}$ is stable in the plot over different but relatively small $k$ values. Note that $\log X_{k,n}$ serves as the threshold for the POT method. Other alternative tail index estimators also available  in the literature, but the performance for these alternative estimators crucially depend on the interplay between the tail index and the so-called second order parameter, and the Hill estimator seems to be most popular due to its simplicity. It is known that (e.g., \cite{Drees+De-Ha+etal:00a}) the Hill plot works most effectively when the underlying loss distribution $F(x)$ is Pareto with 
\begin{equation}
\label{Fbar.pareto}
F(x)=1-\Big{(} \frac{x}{\sigma}\Big{)}^{-\alpha}, \quad x>\sigma>0.
\end{equation} In fact, the Hill estimator implicitly uses the fact that the Pareto variable transforms to a shifted exponential via taking logarithm, and can be seen as an approximate MLE for the exponential distribution with rate $\alpha$, which is the sample mean. Indeed the Hill estimator has minimum mean squared error when the second order parameter governing the rate of convergence of (\ref{regular.varying.ft}) is zero which corresponds to the Pareto case. However, as $F(x)$ departs from Pareto, or equivalently, as the log-transformed data departs from exponential, its performance gets poorer and the plot becomes harder to decipher. Considering that the log-transformed variable of GPD is that of exGPD, not of exponential, the detrimental performance of the Hill plot in identifying the GPD tail is not surprising. \\

To this extent, we propose to use the exGPD directly, which is the correct distribution of the GPD after log-transform. Among other properties of  the exGPD, we employ the second moment result given in (\ref{var.exGPD.trigamma}). The rationale is as follows. When $X$ is a random variable with a heavy tailed distribution, and its GPD realm  starts from some large threshold $u>0$, we assert 
that $X-u|X>u \sim GPD(\sigma, \xi)$ in the spirit of the POT framework. Furthermore, according to the stability property of the GPD,  it is also true that $X-u'|X>u'$ is $GPD(\sigma+\xi(u'-u), \xi)$ distributed, for any higher threshold $u'>u$. Hence, from the definition of the exGPD, we can say that $\log (X-u|X>u) \sim exGPD(\sigma, \xi)$ and $\log (X-u'|X>u') \sim exGPD(\sigma+\xi(u'-u), \xi)$ for $u'>u$. Since our target parameter $\xi$ is invariant under the threshold change within the GPD realm, and the variance of $\log (X-u'|X>u')$ involves $\xi$ only for any $u'>u$ as seen in (\ref{var.exGPD.trigamma}), we are able to estimate $\xi$ by matching the sample variance of log-exceedances to the theoretical one over different thresholds, and produce a series of different $\hat{\xi}$ and plot them. This is always possible because the has a finite variance. A stable area is then identified in the plot, just like in the Hill plot. In particular, we adopt the following algorithm for a given heavy tailed dataset ordered as $X_{n,n} \le ... \le X_{1,n}$ to draw a plot to estimate $\xi$ (or $\alpha=1/\xi$):
\begin{enumerate}
  \item Let $u_i= X_{i,n}$ be the $i$th threshold. 
  \item For all $X_{j,n}$ greater than $u_i$, obtain the exceedances $X_{j,n}-u_i$ and get $\log (X_{j,n}-u_i)$, $j=1,\ldots,i-1$. Equate the sample variance of these $i-1$ values to the theoretical variance in (\ref{var.exGPD.trigamma}) to get an estimate of $\xi$. Denote this by $\tilde{\xi}_{i}$.
  \item Repeat above two steps for $i=n, \dots, 3$ to yield $\tilde{\xi}_{n},\tilde{\xi}_{n-1},\dots, \tilde{\xi}_{3}$. (Note that $i$ stops at 3 because the sample size needs to be at least two for the sample variance to exist)
  \item Plot $(k,  \hat{\xi}_{k} )$ where $\hat{\xi}_{k}= (k-2)^{-1}\sum_{i=3}^k \tilde{\xi}_{i}$, $k=n, \dots, 3.$
\end{enumerate}
We shall call the resulting plot the log variance (LV) plot as it is based on the sample variance of the log-transformed exceedances. Note that Step 4 uses an average of $\tilde{\xi}$ value to smooth volatile sample variance values in the plot. This is similar to the idea found in \cite{Resni+Stari:97a} where the averaged Hill estimators are plotted to reduce the variability of the standard plot.  In the following two subsections we compare the performances of the proposed algorithm and the Hill estimator for simulated and real datasets from heavy tailed distributions, respectively.

\subsection{Comparison between Hill plot and LV plotSimulated data from GPD and GEV}
\begin{figure}[H]
\centering
\includegraphics[scale=0.23]
{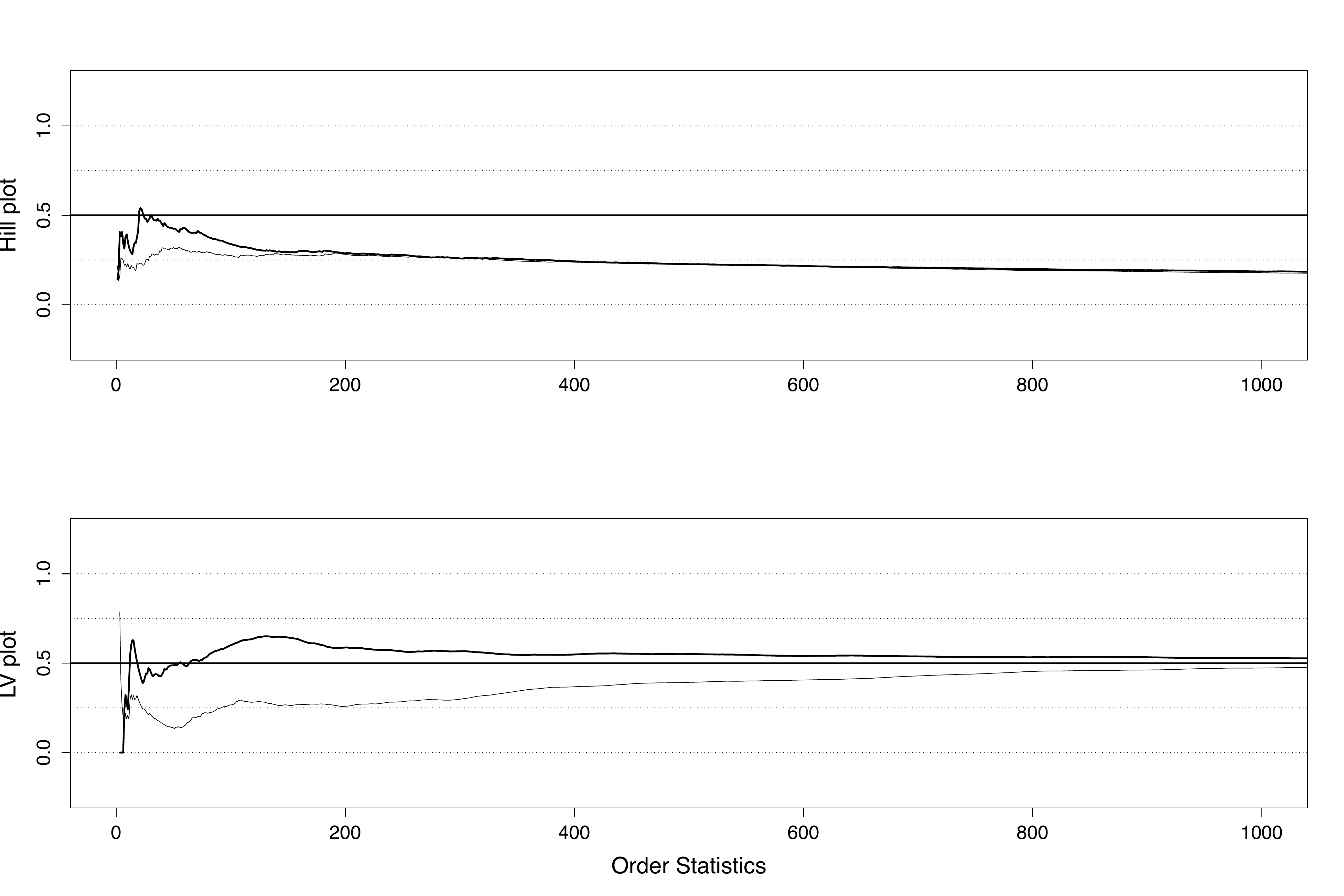}
\caption{Comparison of $\xi$ for GPD with $\mu=10$, $\sigma=1$, $\xi = 0.5$}
\label{fig:mu10xi05}
\end{figure}

\begin{figure}[H]
\centering
\includegraphics[scale=0.23]
{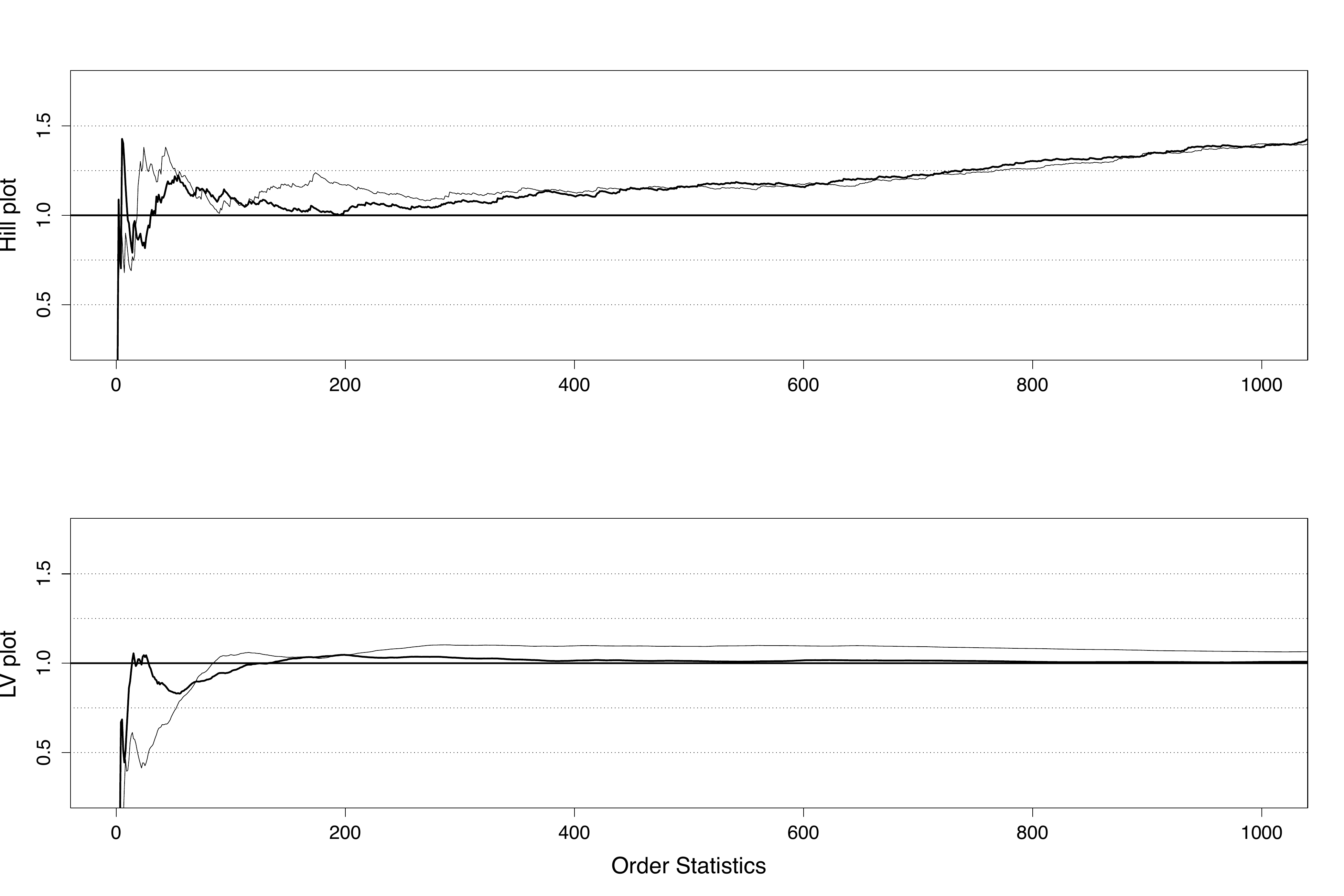}
\caption{Comparison of $\xi$ for GPD with $\mu=0$, $\sigma=1$, $\xi = 1$}
\label{fig:mu0xi1}
\end{figure}

\begin{figure}[H]
\centering
\includegraphics[scale=0.23]
{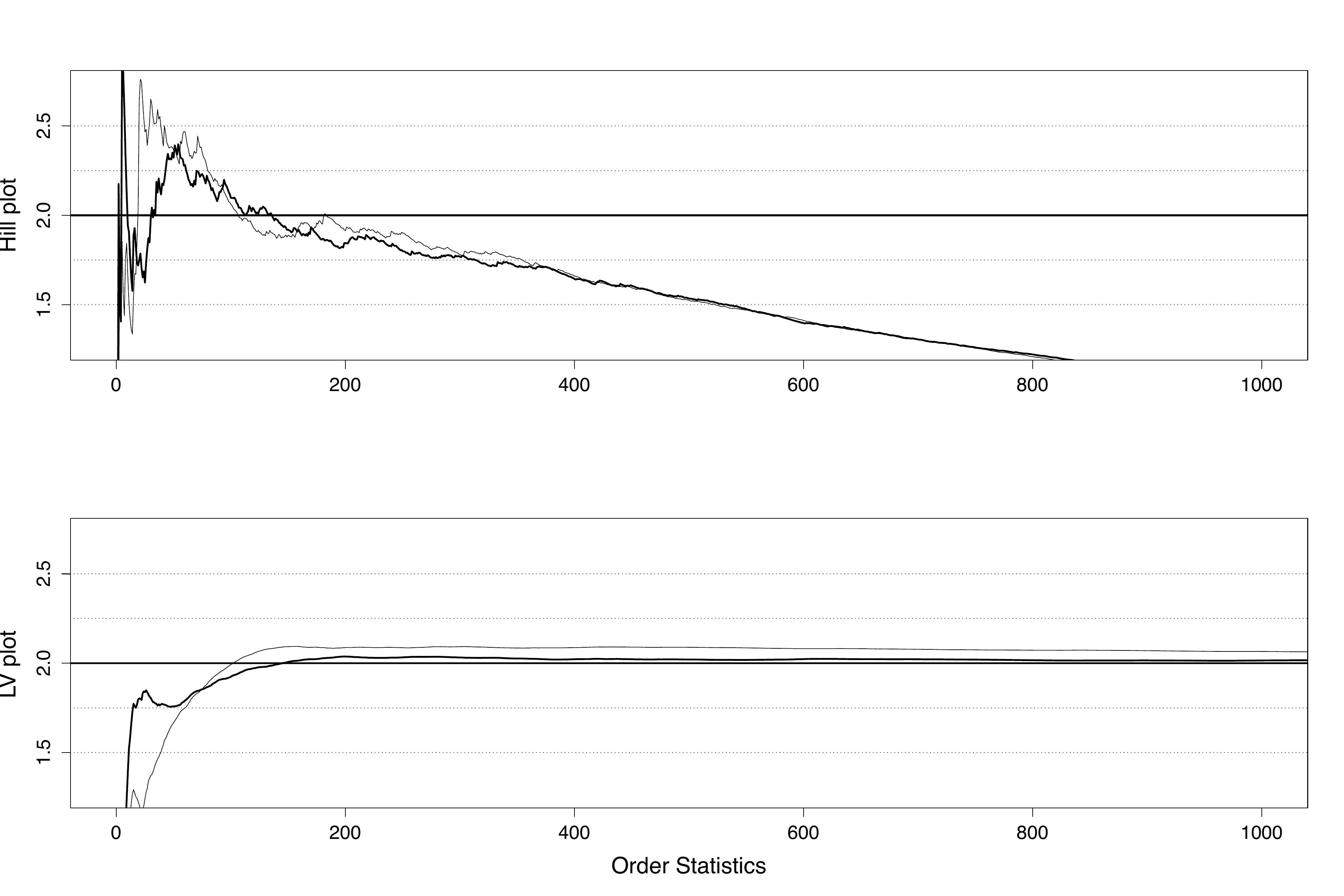}
\caption{Comparison of $\xi$ for GPD with $\mu=10$, $\sigma=1$, $\xi = 2$}
\label{fig:mu10xi2}
\end{figure}

\begin{figure}[H]
\centering
\includegraphics[scale=0.23]
{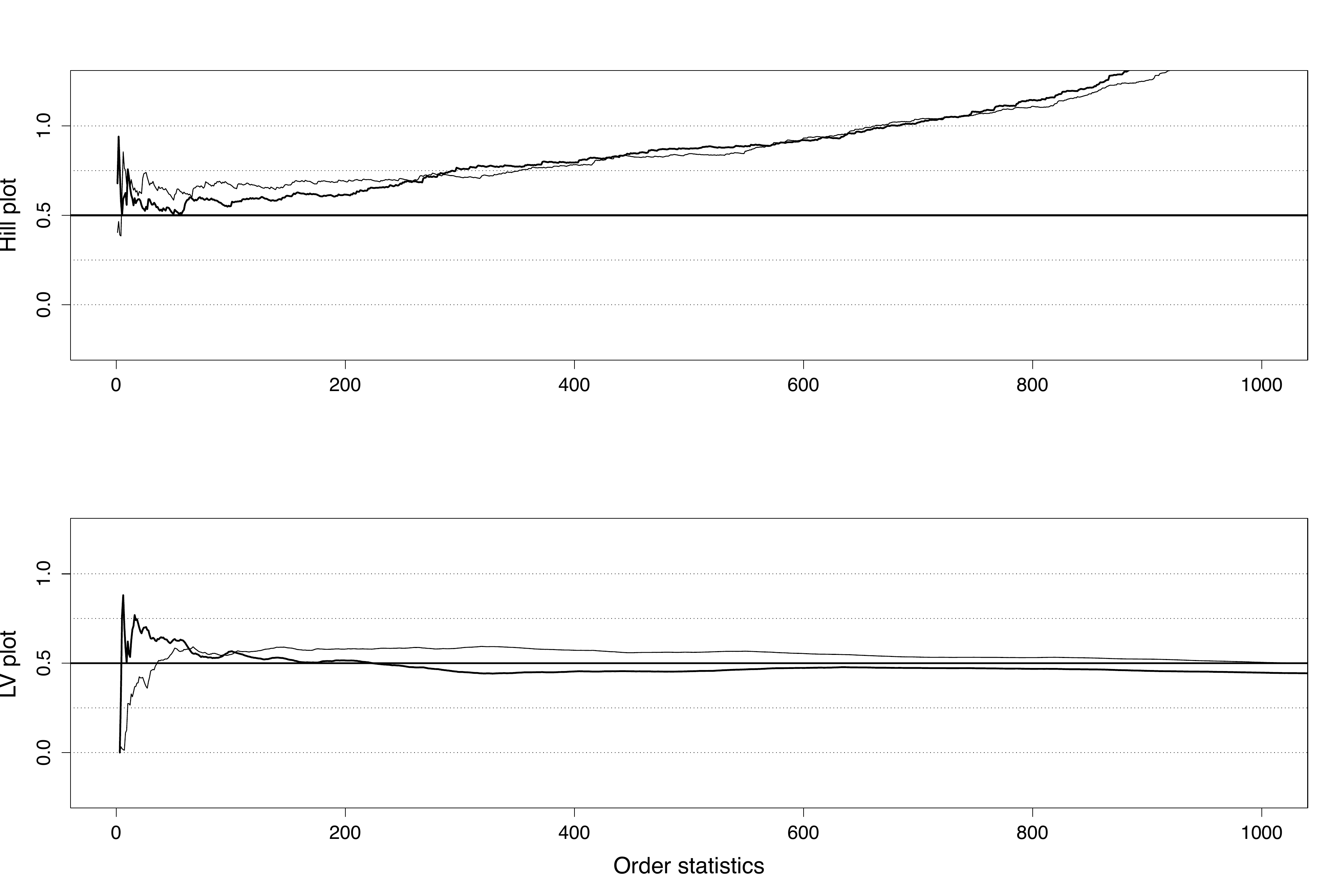}
\caption{Comparison of $\xi$ for GEV with $\mu=0$, $\sigma=1$, $\xi = 0.5$}
\label{fig:mu0xi05gev}
\end{figure}

\begin{figure}[H]
\centering
\includegraphics[scale=0.23]
{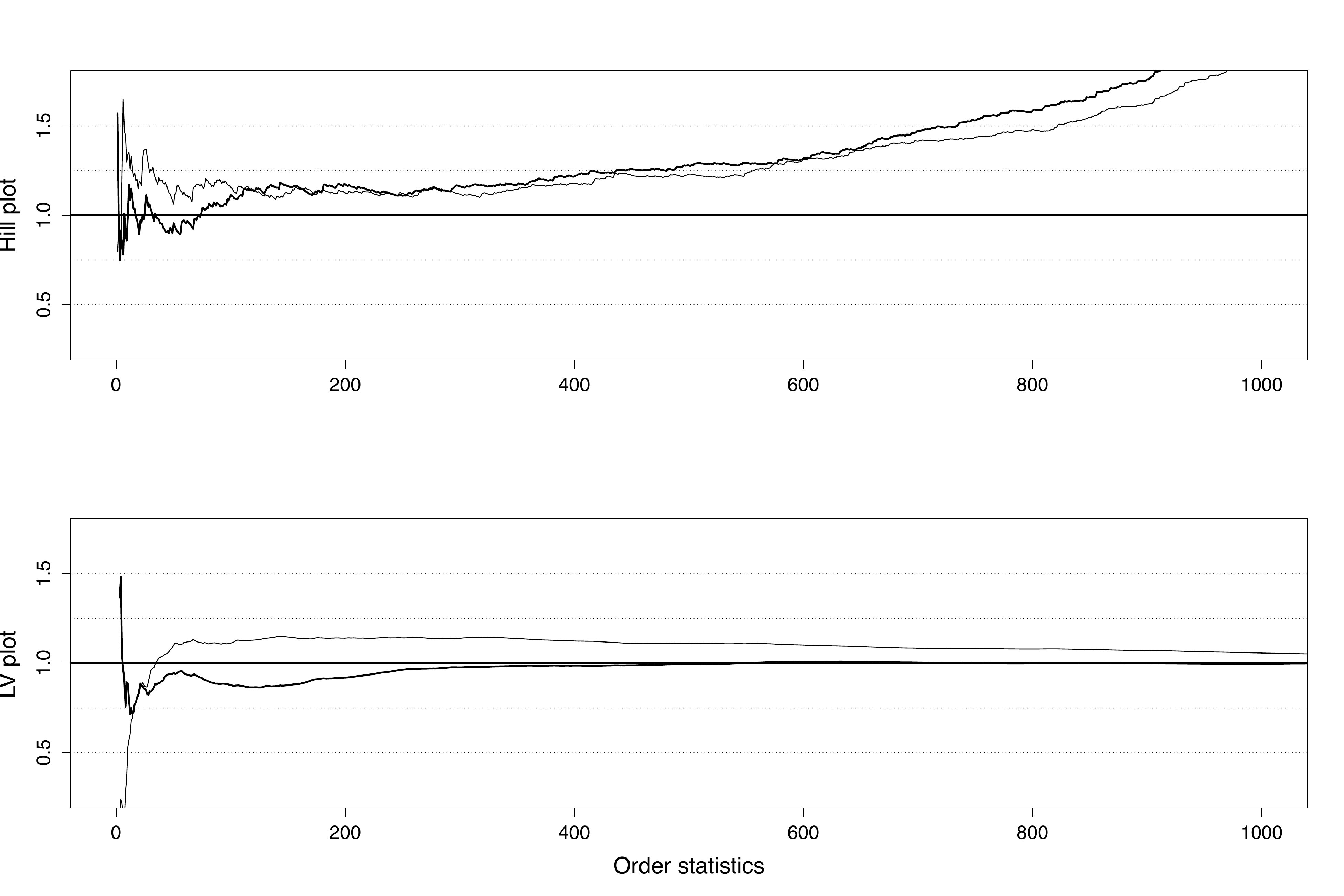}
\caption{Comparison of $\xi$ for GEV with $\mu=0$, $\sigma=1$, $\xi = 1$}
\label{fig:mu0xi1gev}
\end{figure}

\begin{figure}[H]
\centering
\includegraphics[scale=0.23]
{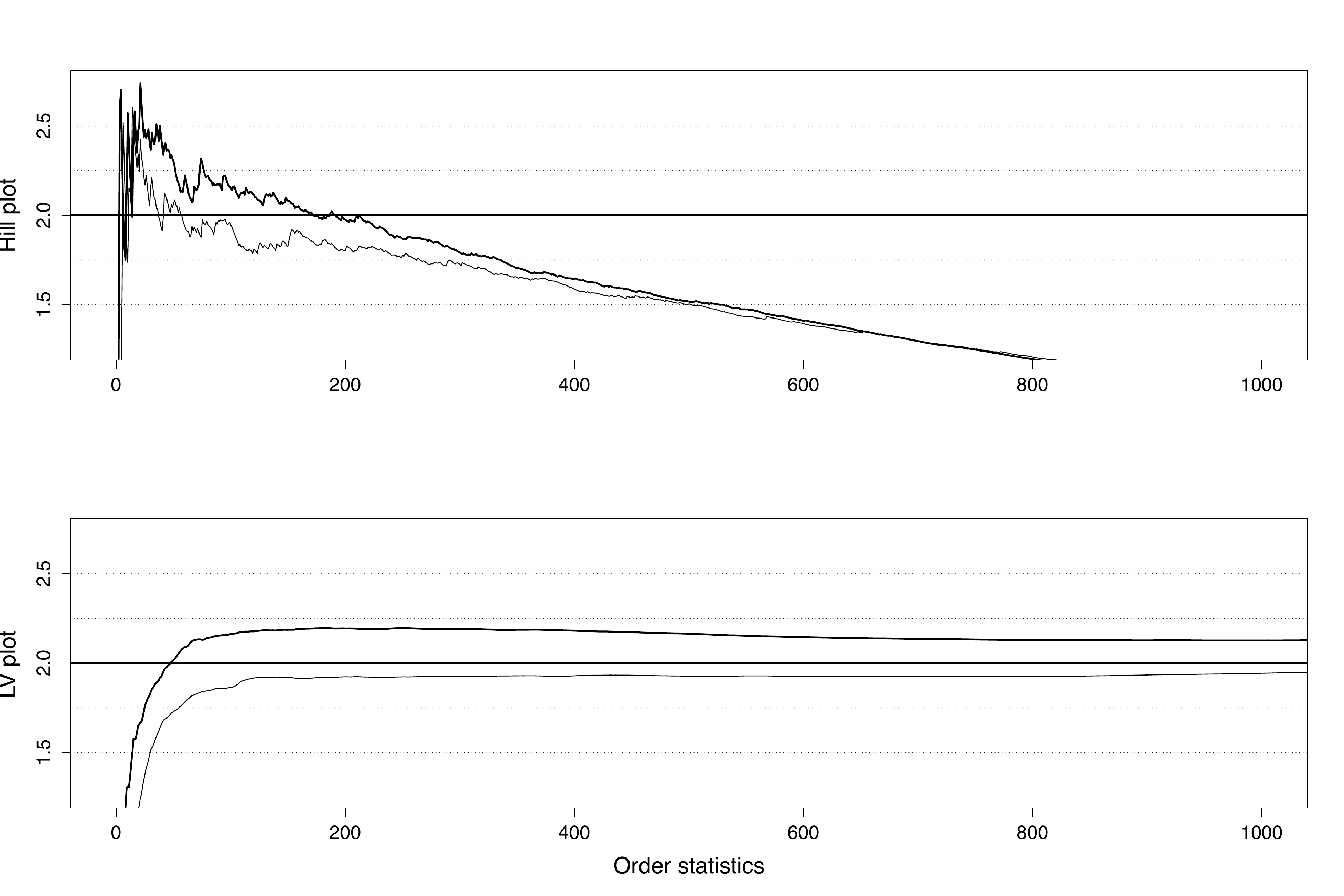}
\caption{Comparison of $\xi$ for GEV with $\mu=10$, $\sigma=1$, $\xi = 2$}
\label{fig:mu0xi2gev}
\end{figure}

Our simulated datasets are from the GPD and GEV, two most prominent distributions in the EVT studies. To compare the performance of the Hill and LV plots, we consider several GPD and GEV distributions with different choices for parameters $\xi, \sigma$ and $\mu$. The location parameter $\mu$ has also been added in our simulation study to allow diverse heavy tail shapes. As our LV plot algorithm above is invariant under the location shifting in $X$, it can handle these simulated data with no modification. For each GPD and GEV chosen, we generate samples of size 2,000 and present two selected sample paths for illustration.  We have tried larger sample sizes and more repetitions, but the conclusions are similar. As distributions with finite upper bounds hardly concern us in real applications, we focus on $\xi>0$ cases only and negative estimates of $\xi$ are accordingly replaced with 0 in computing $\tilde{\xi}$ in the algorithm; this is also consistent with the Hill plot which works only for $\xi>0$. The results are depicted in Figures \ref{fig:mu10xi05} to \ref{fig:mu0xi2gev} where the first three figures used GPD samples and the last three used GEV samples. In each figure, the upper and lower panels represent the Hill and LV plots, respectively. The vertical line shows the value of estimated $\xi$. The horizontal axis stands for the ordered data count from the largest observation; for example, 200 in the horizontal axis means the largest 200 observations in calculating the estimate. Note that we have limited the range of the figures to the largest 1,000 observations even though the same size is 2,000.  So, basically the left end area in each figure is the region where the GPD phenomenon takes place and we wish to find the tail index in this area.  \\

From these six figures we see that, depending the combination of the parameter sets, the Hill plot can be stable but also can have a trend of increasing or decreasing over the whole dataset, confirming that the Hill plot should be read only over  a small fraction of the dataset on the left end side of the figure, or the tail region. However, since we do not know exactly where the tail starts, any trend makes the Hill plot harder to decipher. Practical experience (e.g., Section 7.2 in \cite{mcne+frey+embr:05}) suggests the 5\% upper tail to be used to determine the tail index with the Hill plot. When focused on this region in each figure, the Hill plot is sometimes able to reveal the true tail index (e.g., Figures  \ref{fig:mu0xi1} and   \ref{fig:mu0xi05gev}), but in other cases there is hardly any stable region in the plot (e.g., Figures \ref{fig:mu10xi2}  and \ref{fig:mu0xi2gev}). Thus the Hill plot, despite its popularity, remains difficult to use in practical situations. 
If we turn to the varaince plot, finding the tail index is relatively easier because it has been smoothed via averaging, though both plots can be qualitatively poor when the sample itself has large variation in the tail. For example, in Figures  \ref{fig:mu10xi05} and \ref{fig:mu0xi1gev}, both plots are relatively stable at the left end region but they point at the wrong tail index. In Figures \ref{fig:mu10xi2} and \ref{fig:mu0xi2gev}, both plots are either unstable or biased at the left end area; in particular, the Hill plot is highly volatile and the LV plot has a considerable bias. Compared to the Hill plot, however, a marked advantage of the LV plot throughout all figures is that it exhibits a quite stable region in not so extreme quantile levels without any substantial trend, and in this region the estimated value is fairly close to the true tail index, again, relative to the Hill plot. These observations suggest that the LV plot can be most useful when one reads the plot between, say, upper 5\% and 20\% quantile regions. This is advantageous for practitioners to interpret the plot as the volatility of extreme quantiles can be avoided in this non-extreme quantile range.

\subsection{Real datasets}
We now repeat the same exercise for two actual datasets publicaly available: the Danish fire data and the BMW stock price log-returns. The famous Danish fire data, collected at Copenhagen Reinsurance, comprise 2,167 fire losses between 1980 and 1990, both years inclusive. The numbers have been adjusted for inflation to reflect 1985 values and are expressed in millions of Danish Krone. The BMW data consists of daily  log-returns during January 2, 1973 and July 23, 1996, with a total of 6,146 observations. Both datasets can be obtained from \verb"fExtremes" package in R software, and further description of these datasets can be found in \cite{EmbrechtsModellingExtremalEvents}.\\

\begin{figure}[H]
\centering
\includegraphics[scale=0.23]
{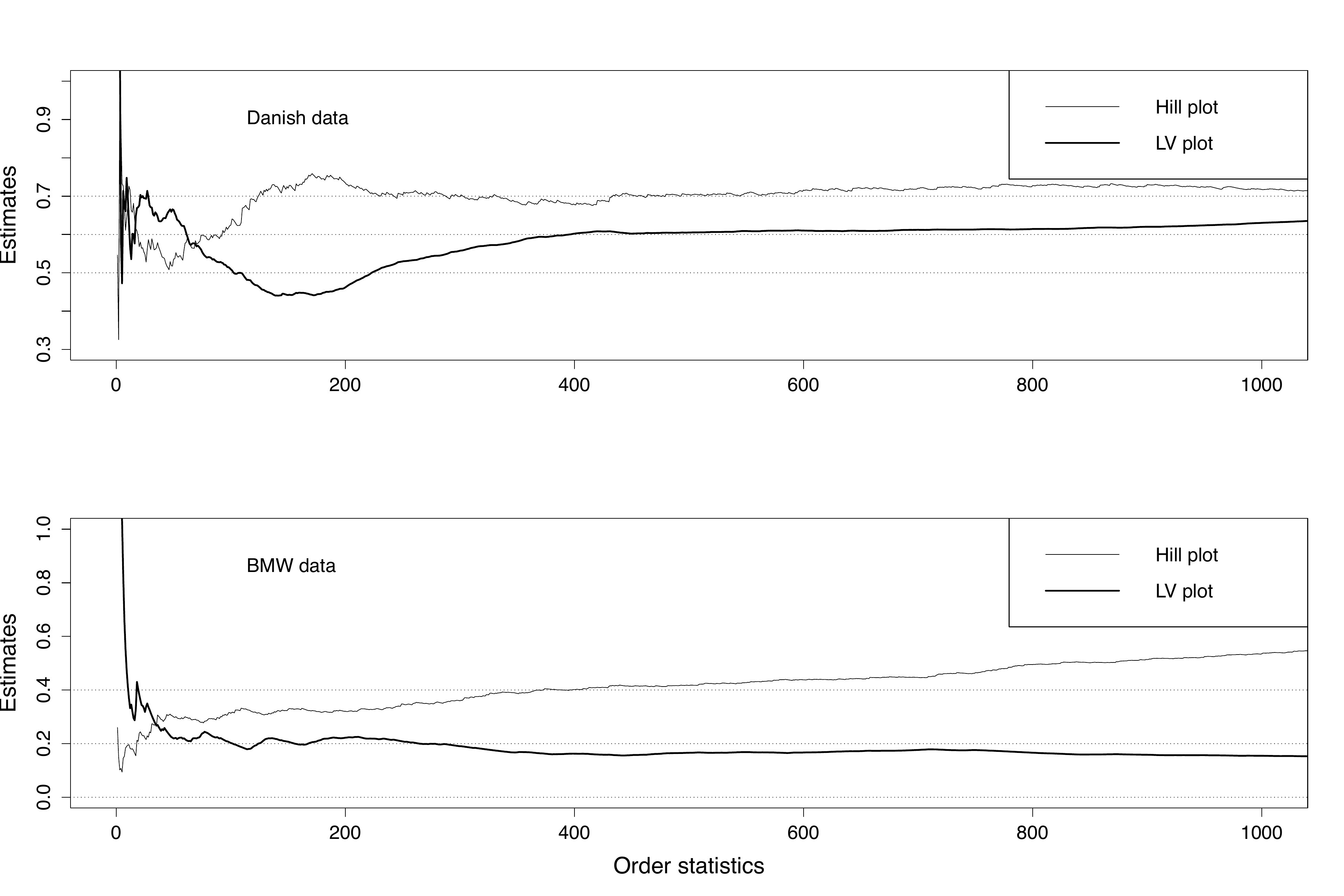}
\caption{Plot comparison of $\xi$ for Danish (top) and BMW data (bottom)}
\label{fig:danish_bmw}
\end{figure}

The Hill and LV plots for these two datasets are presented in Figure \ref{fig:danish_bmw}. For the Danish data (top panel), previously two choices for the tail index have been suggested. According to Section 6.4 and 6.5 of \cite{EmbrechtsModellingExtremalEvents}, the first threshold is at a loss of $10$ (109 exceedances) with the corresponding Hill estimate $\hat{\alpha}^{-1}=0.618$, indicating that only the first moment exists. For the second choice, the GPD threshold is set at $18$ (47 exceedances) and  the Hill estimate is $\hat{\alpha}^{-1}=0.497$, suggesting the existence of the first two moments. So essentially ${\alpha}^{-1}={\xi}$ is indicated to be in $[0.5, 0.62]$ by the Hill plot. However, the chosen GPD thresholds, 10 and 18, were suggested based on the mean excess plot because the Hill plot itself does not show any stable region within the upper 5\% quantile range as seen from the figure. When we use the LV plot, the path is more tamed and easier to investigate. In fact, following the recommendation from the previous simulation study, we search the plot between upper 5\% and 20\% quantile regions and estimate that the true $\xi$ lies in $[0.45, 0.65]$. This is largely in agreement with the previously suggested interval above, though our interval is slightly wider.  Regarding the BMW data, we first observe that both plots are relatively stable at the left end, implying that the extreme tail variability is not substantial. The Hill plot suggests that $\xi$ lies roughly in $[0.25, 0.3]$ by inspecting the upper 5\% quantile range, indicating that only the first 3 or 4 moments are existent. In contrast, the LV plot indicates a slightly different interval of $[0.19, 0.25]$ based on the overall plot, indicating a less heavier tail with finite moments up to 4 or 5. As the LV plot in this data has a very narrow interval for $\xi$ value over the entire data range, we have more confidence in estimating the tail index. We comment that, as observed from previous simulation studies,  the interval for $\xi$ in the LV plot can be identified in a much clearer manner for these real datasets  compared to the Hill plot.  \\

\section{Concluding remarks}
In this paper we propose  and study the Exponentiated Generalized Pareto Distribution (exGPD), which is created via log-transform of the Generalized Pareto Distribution (GPD), an influential distribution in  Extreme Value Theory. For this distribution  we derive various distributional quantities, including the moment generating function, tail risk measures and quantities related to order statistics. As an application we develop a  plot as an alternative to the Hill plot to identify the tail index of heavy tailed datasets. The proposed plot is based on the idea of the sample variance of log exceedances to be matched to the variance of the exGPD. Through various numerical illustrations with both simulated and actual datasets, it is shown that the proposed plot works reasonably well compared to the Hill plot, elucidating the usefulness of the exGPD.

\end{document}